\documentclass{article}

\usepackage{microtype}
\usepackage{graphicx}
\usepackage{subfigure}
\usepackage{booktabs} 

\usepackage{hyperref}

\usepackage[accepted]{my_icml2024}

\usepackage{amsmath}
\usepackage{amssymb}
\usepackage{mathtools}
\usepackage{amsthm}

\usepackage{thmtools,thm-restate}
\usepackage{mdframed}

\usepackage[capitalize,noabbrev]{cleveref}

\newcommand{\mat}[1]{\mathbf{#1}}
\def\argmin{\mathop{\rm argmin}}

\def\argmax{\mathop{\rm argmax}}

\def\rank{{\rm rank \,}}
\def\smax{{\rm smax \,}}

\def\Diag{{\rm Diag \,}}
\def\Def{\stackrel{\mathrm{def}}{=}}
\def\im{{\rm Im \,}}
\def\ker{{\rm Ker \,}}
\def\inter{{\rm int \,}}

\def\dom{{\rm dom \,}}

\def\R{\mathbb{R}}

\def\F{\mathcal{F}}

\providecommand{\0}{\boldsymbol{0}}

\renewcommand{\aa}{\boldsymbol{a}}
\providecommand{\bb}{\boldsymbol{b}}
\providecommand{\cc}{\boldsymbol{c}}

\let\ggg\gg
\renewcommand{\gg}{\boldsymbol{g}}
\providecommand{\hh}{\boldsymbol{h}}

\renewcommand{\ss}{\boldsymbol{s}}

\providecommand{\uu}{\boldsymbol{u}}
\providecommand{\vv}{\boldsymbol{v}}

\providecommand{\xx}{\boldsymbol{x}}
\providecommand{\yy}{\boldsymbol{y}}

\providecommand{\cO}{\mathcal{O}}

\providecommand{\cO}{\mathcal{O}}

\def\beq{\begin{equation}}
\def\eeq{\end{equation}}
\def\ba{\begin{array}}
	\def\ea{\end{array}}
\def\la{\langle}
\def\ra{\rangle}

\def\BT{\begin{theorem}}
	\def\ET{\end{theorem}}
\def\BL{\begin{lemma}}
	\def\EL{\end{lemma}}
\def\BC{\begin{corollary}}
	\def\EC{\end{corollary}}
\def\BE{\begin{example}}
	\def\EE{\end{example}}
\def\BD{\begin{definition}}
	\def\ED{\end{definition}}
\def\BR{\begin{remark}}
	\def\ER{\end{remark}}
\def\BAS{\begin{assumption}}
	\def\EAS{\end{assumption}}
\def\BP{\begin{proposition}}
	\def\EP{\end{proposition}}

\def\BTF{\begin{framedtheorem}}
	\def\ETF{\end{framedtheorem}}

\definecolor{mydarkgreen}{RGB}{39,130,67}

\definecolor{mydarkred}{RGB}{192,47,25}

\theoremstyle{plain}
\newtheorem{theorem}{Theorem}[section]
\newtheorem{proposition}[theorem]{Proposition}
\newtheorem{lemma}[theorem]{Lemma}
\newtheorem{corollary}[theorem]{Corollary}
\theoremstyle{definition}
\newtheorem{definition}[theorem]{Definition}
\newtheorem{assumption}[theorem]{Assumption}
\theoremstyle{remark}
\newtheorem{remark}[theorem]{Remark}
\theoremstyle{example}
\newtheorem{example}[theorem]{Example}

\newmdtheoremenv{framedtheorem}[theorem]{Theorem}
\newmdtheoremenv{framedlemma}[theorem]{Lemma}
\newmdtheoremenv{framedexample}[theorem]{Example}
\newmdtheoremenv{framedassumption}[theorem]{Assumption}
\newmdtheoremenv{framedproposition}[theorem]{Proposition}

\usepackage[textsize=tiny]{todonotes}

\icmltitlerunning{Spectral Preconditioning
	for Gradient Methods on Graded Non-convex Functions}

\begin{document}
	
	\twocolumn[
	\icmltitle{Spectral Preconditioning
		for Gradient Methods\\on Graded Non-convex Functions}

	
	
	
	\begin{icmlauthorlist}
		\icmlauthor{Nikita Doikov}{EPFL}
		\icmlauthor{Sebastian U. Stich}{CISPA}
		\icmlauthor{Martin Jaggi}{EPFL}
	\end{icmlauthorlist}
	
	\icmlaffiliation{EPFL}{Machine Learning and Optimization Laboratory (MLO), EPFL, Lausanne, Switzerland}
	\icmlaffiliation{CISPA}{CISPA Helmholtz Center for Information Security, Saarbrücken, Germany}
	
	\icmlcorrespondingauthor{Nikita Doikov}{nikita.doikov@epfl.ch}
	\icmlcorrespondingauthor{Sebastian U. Stich}{stich@cispa.de}
	\icmlcorrespondingauthor{Martin Jaggi}{martin.jaggi@epfl.ch}
	
	\icmlkeywords{Gradient Methods, Preconditioning, Hessian Spectrum, Global Convergence, Non-convex Optimization}
	
	\vskip 0.3in
	]
	
	
	
	\printAffiliationsAndNotice{}  
	
	\begin{abstract}
        
        The performance of optimization methods is often tied to the spectrum of the objective Hessian.
        Yet, conventional assumptions, such as smoothness, do often not enable us to make finely-grained convergence statements---particularly not for non-convex problems. 
        Striving for a more intricate characterization of complexity, we introduce a unique concept termed \emph{graded non-convexity}. 
        This allows to partition the class of non-convex problems into a nested chain of subclasses. 
        Interestingly, many traditional non-convex objectives, including partially convex problems, matrix factorizations, and neural networks, fall within these subclasses. 
        As a second contribution, we propose gradient methods with spectral preconditioning, which employ inexact top eigenvectors of the Hessian to address the ill-conditioning of the problem, contingent on the grade.
        Our analysis reveals that these new methods provide provably superior convergence rates compared to basic gradient descent on applicable problem classes, particularly when large gaps exist between the top eigenvalues of the Hessian. 
        Our theory is validated by numerical experiments executed on multiple practical machine learning problems.

	\end{abstract}

	\section{Introduction}
	\label{SectionIntroduction}

	\paragraph{Motivation.}

		\begin{figure}[h!]
		\begin{center}
			\includegraphics[scale=0.40]{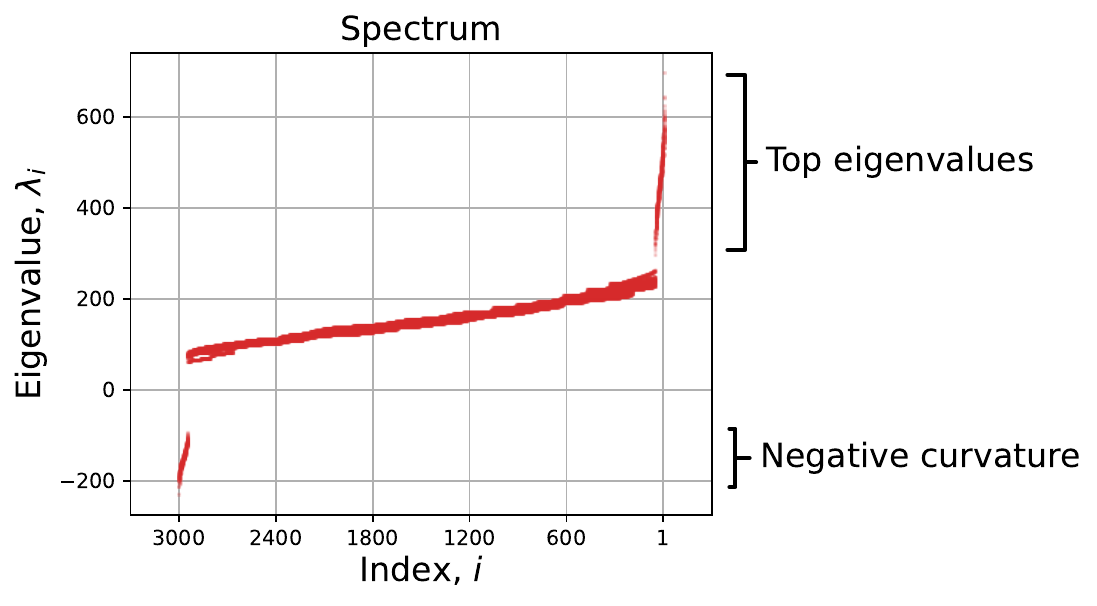}
		\end{center}\vspace{-1em}
		\caption{ \small Spectrum of the Hessians for the matrix factorization problem (Example~\ref{ExampleMF}) 
			of the optimization dimension $n = 3000$,
			for $10$ random objectives.  }
		\label{FigSpectrum}
	\end{figure}

	The gradient method is an important and attractive tool 
	for solving large-scale optimization problems.
	It has very cheap cost of every iteration
	and well established convergence guarantees,
	that hold starting from an arbitrary initial point
	and for a wide family of problem classes,
	including convex and non-convex problems.
	However, the major drawback of the gradient method
	remains to be its slow rate of convergence:
	for solving modern optimization problems
	up to a reasonable accuracy level,
	it is often required to do a lot of gradient steps, due 
	to \textit{ill-conditioning} of the problem.
	
	In order to improve the gradient direction,
	we can multiply the gradient by  a specifically crafted matrix called
	\textit{preconditioner}, which should adjust the method to the right geometry
	of the problem. However, finding a good preconditioning
	with strong theoretical guarantees
	is not easy, especially for non-convex problems. In this work, 
	we propose a new family of \textit{spectral preconditioners}
	that rely on an additional refined information about the function class.
	As a by-product, we establish convergence rates, that are provably better than those
	of the gradient methods.

	\begin{table*}[h!]
		\centering
		\small
		\renewcommand{\arraystretch}{1.6}
		\begin{tabular}{r|c|c|c|c}
			\toprule
			Algorithm & Preconditioning & Non-convex Complexity & Strongly Convex & Arithmetic Cost \\
			\hline
			The Gradient Method  & $\mat{H} := \mat{I}$, \; $\tau = 0$  & 
			$ \cO\bigl( 
			{\color{mydarkred} \boldsymbol{\frac{\lambda_1}{\varepsilon^{2}} }} \bigr) \cdot F_0 $ & 
			$\tilde{\cO}\bigl(  
			{\color{mydarkred} \boldsymbol{ \frac{\lambda_1}{\lambda_n} }}  \bigr)$ & 
			$\cO\bigl( {\color{mydarkgreen} \boldsymbol{n}  } \bigr)$ \\
			\hline
			Spectral Preconditioning \textbf{(ours)} &
			$\mat{H} \approx \nabla^2_{\tau} f$, \; $1 \leq \tau \leq n$ &
			$\cO\bigl(  \,
					{\color{mydarkgreen} \boldsymbol{
					\frac{  \lambda_{\tau + 1}  }{\varepsilon^{2}}  }}
				+  \frac{L^{1/2}}{\varepsilon^{3/2}} 
			
			\, \bigr)
			\cdot F_0$ & 
			$\tilde{\cO}\bigl(  
			{\color{mydarkgreen} \boldsymbol{
			\frac{  \lambda_{\tau + 1}  }{\lambda_n} }} + MD   \bigr)$ &
			$\cO\bigl( 
			{\color{mydarkgreen} \boldsymbol{ \tau^2 n + \tau^3 } }  \bigr)$
			\\
			\hline
			Newton's Method & $\mat{H} := \nabla^2 f$, \; $\tau = n$ &
			$\cO\bigl(  \frac{L^{1/2}}{\varepsilon^{3/2}} \bigr) \cdot F_0$ &
			$\tilde{\cO}\bigl( MD \bigr)$ &
			$\cO\bigl( 
				{\color{mydarkred} \boldsymbol{  n^3  }} 
				\bigr)$ \\
			\bottomrule
		\end{tabular}
		\caption{Global complexity bounds for our method as compared to the classic gradient method
		and the regularized Newton methods. We denote by $\lambda_i$ a uniform bound
		for the $i$th eigenvalue of the Hessian, sorted in a non-ascending order:  $ \lambda_1 \geq \ldots \geq \lambda_n$.
		We denote by $L$ the Lipschitz constant of the Hessian and, for convex objectives, $M$ is the constant of quasi-self-concordance
		(see Section~\ref{SectionConvex}).
	 	$F_0 := f(\xx_0) - f^{\star}$, and $D$ is the diameter of the initial sublevel set. 
 		We see that using the spectral preconditioner of order $\tau$, we cut the top $\tau$ eigenvalues of the spectrum,
 		which is the most significant complexity factor.
 		We present the state-of-the-art global complexities from \cite{nesterov2006cubic,doikov2023minimizing}
 		for Newton's Method with cubic and gradient regularizations.
 		}
		\label{TableComplexities}
	\end{table*}

	Optimization theory suggests that the main 
	complexity parameters that affect the rate 
	of convergence for gradient methods are the spectrum of the Hessian $\nabla^2 f$
	and its extremal characteristics, such as the bound
	for the \textit{maximal eigenvalue} (the Lipschitz constant of the gradient)
	or the \textit{condition number} (the ratio of the largest and smallest eigenvalues)
	\cite{nemirovskii1983problem,nesterov2018lectures}.
	Moreover, some of the fundamental properties
	of the objective function, such as \textit{convexity},
	\textit{weak convexity},
	or \textit{strong convexity},
	that distinguish between all optimization problems
	and globally solvable ones, can be defined in terms
	of lower bounds on the spectrum. Thus (for a twice differentiable
	function),\vspace{-2mm}
	\beq \label{ConvexDef}
	\ba{rcl}
	f \; \text{is convex} \quad \Leftrightarrow \quad
	\nabla^2 f \; \succeq \; \0.
	\ea
	\eeq
	At the same time, for problems with specific structure,
	the worst-case guarantees for convergence of gradient methods obtained considering
	only the the largest and smallest eigenvalue of the Hessian
	can be too pessimistic.
	Indeed, we see that in practice,
	the distribution of eigenvalues can be quite specific,
	with a relatively small amount of top eigenvalues
	that are much larger than the others (Fig.~\ref{FigSpectrum}).
	Consequently, any a priori information on the structure of the Hessian
	can be significant from the optimization perspective. 
	Ultimately, we want to have algorithms that are able to benefit from 
	this knowledge, achieving faster rates
	when the distribution of eigenvalues is far from uniform.

	\begin{figure}[h!]
		\begin{center}
			\includegraphics[scale=0.40]{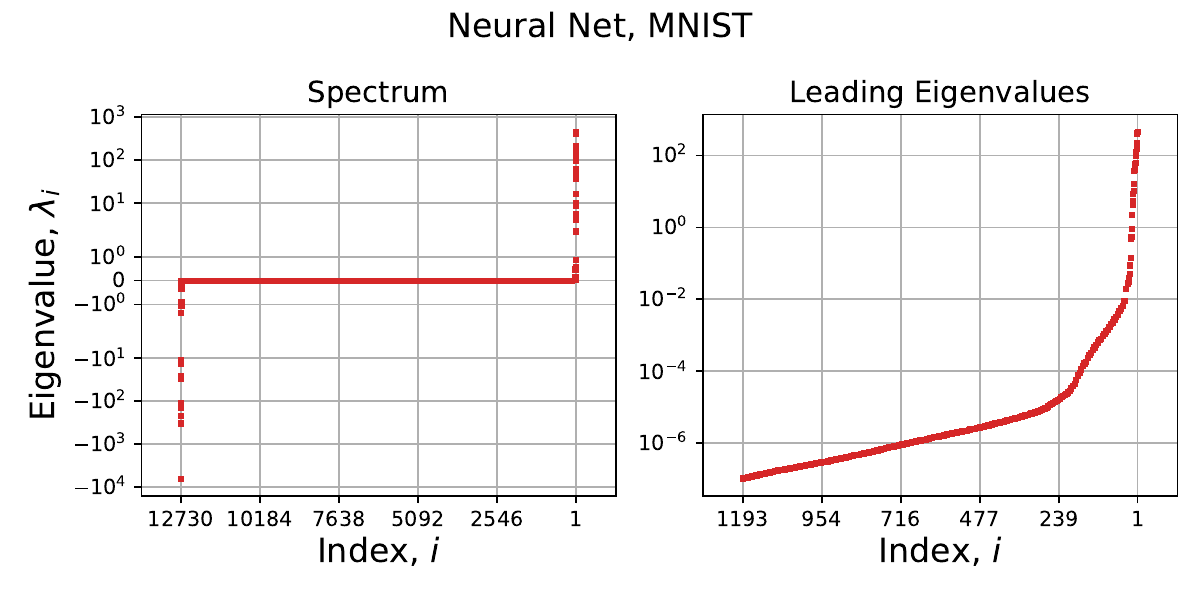}
		\end{center}\vspace{-1em}
		\caption{ \small Left: spectrum of the Hessians for a 
			two-layer fully connected neural network
			trained on MNIST dataset. Right: zoomed top eigenvalues.
			The total number of parameters is $n = 12730$
			(the first layer: $12560$, the second layer: $170$).
			We see that the dimension $\tau$ of the subspace with positive eigenvalues
			is much bigger than the dimension of the last layer. However, there are only a few
			eigenvalues that are significantly larger than the others.
		}
		\label{FigSpectrumMNIST}
	\end{figure}

	\paragraph{Contributions.}

	In this work, we develop spectral preconditioning
	for the gradient methods that is able to tackle
	non-convex problems with \textit{highly non-uniform} and clustered spectrum
	of the Hessian, that we often observe in practice.
	For that, we propose to make a step back from the common
	dichotomy between  convex and non-convex problems.
	We introduce a new notion of \textit{graded non-convex}
	functions,
	which granulate the class of all non-convex problems
	into nested family of subclasses. We say, for some integer $\tau \geq 1$
	\beq \label{GradedDef}
	\ba{rcl}
	f \; \text{is non-convex of grade} \; \tau \quad \Leftrightarrow \quad
	\nabla^2_{\tau} f \; \succeq \; \0,
	\ea
	\eeq
	where $\nabla^2_{\tau} f$ is composed by the top $\tau$ eigenvectors 
	of the full Hessian (see \eqref{HessOrder} for the formal definition).
	For $\tau = n$, where $n$ is dimension of the problem,
	we obtain the entire Hessian $\nabla^2 f$ and
	inequality~\eqref{GradedDef} means the standard convexity~\eqref{ConvexDef}.
	It appears that for many practical non-convex problems,
	this condition is satisfied at least for some small $\tau$.
	For example, for deep neural networks with convex loss it is satisfied 
	at least
	for $\tau \geq d$, where $d$ is the dimension of the last layer
	(however, in practice we observe that the actual values of $\tau$ can be much bigger, see Fig.~\ref{FigSpectrumMNIST}).
	
	Inequality~\eqref{GradedDef} provides us
	with a certain \textit{convex subspace} at each point 
	(such that our function is locally convex in this subspace, see Fig.~\ref{FigureSurface}),
	which looks very attractive from the optimization standpoint.
	Furthermore, it contains the leading eigenvalues, which constitute the primary computational burden for first-order methods.
	
	Based on our problem class,
	we propose to use a positive definite matrix $\mat{H} \approx \nabla^2_{\tau} f(\xx)$ 
	as a natural preconditioner for our method,
	that we call \textit{spectral preconditioning}:
	\beq \label{IntroStep}
	\ba{rcl}
	\xx^+ & = & \xx - \bigl( \mat{H} + \alpha \mat{I} \bigr)^{-1} \nabla f(\xx),
	\ea
	\eeq
	where and $\alpha \geq 0$ is a regularization parameter.
	For these iterations, we establish strong convergence guarantees,
	that improve with increasing the parameter~$\tau$ (see Table~\ref{TableComplexities}).
	
	Therefore, the method with any $\tau \geq 1$ works \textit{provably better}
	than the basic gradient method ($\tau = 0$)
	in terms of the rate of convergence.
	For $\tau = n$ (the full Hessian), our algorithm becomes the
	regularized Newton method
	with the best global complexity bounds known in the literature.
	However, the most efficient version of our method
	corresponds to the case of small $\tau \approx 10^0-10^2$, 
	when estimating the matrix~$\mat{H}$
	can be done efficiently with the hot-start \textit{power method}.

	\paragraph{Related Work.}

	Preconditioning is an important tool in numerical analysis
	and optimization \cite{nocedal2006numerical}.
	The basic example is preconditioning of the conjugate gradient method \cite{hestenes1952methods}
	for solving a system of linear equations.
	The choice of the right preconditioner is a difficult task and it often
	depends on an a-priori knowledge on the problem structure,
	as, for example, the Laplacian preconditioning for 
	the graph-induced problems \cite{spielman2004nearly, vaidya1991solving}
	or for the systems involving partial differential equations \cite{mardal2011preconditioning}. 
	
	In optimization, a powerful approach for preconditioning
	that works for general problems is called the Newton method \cite{polyak2007newton}.
	It uses the Hessian matrix $\nabla^2 f$ to alleviate the impact of
	ill-conditioned characteristics of the Hessian spectrum.
	The modern versions of this method provides us with the global rates
	of convergence that are significantly better than those of the first-order algorithms
	\cite{nesterov2006cubic,cartis2011adaptive1,cartis2011adaptive2,grapiglia2017regularized,karimireddy2018global,doikov2022super}.
	However, computing and inverting the Hessian matrix is 
	prohibitively expensive in terms of the arithmetic operations 
	and memory usage to store big matrices.
	
	The common approximate technique for improving the arithmetic cost of the methods with the full Hessian
	is called the \textit{quasi-Newton methods}, such as SR1, DFP, BFGS, L-BFGS and others 
	\cite{nocedal2006numerical}.
	Despite these methods show an outstanding performance 
	on practical problems of a moderate dimension,
	it is a significant challenge to establish
	a rigorous theory of convergence
	for \textit{quasi-second-order methods},
	which would provably benefit globally from inexact Hessian information,
	while the local non-asymptotic theory of convergence
	for the classical quasi-Newton methods has emerged only very
	 recently~\cite{rodomanov2021new,rodomanov2021rates,jin2023non}.
	Note that employing 
	a naive and straightforward line-search or damping
	approach in the basic Newton method can have as slow convergence
	as the plain gradient descent \cite{cartis2013example}.
	See also \cite{hanzely2022damped,kamzolov2023cubic,jiang2023accelerated}
	for the recent advancement in the global analysis of the damped and quasi-Newton methods
	for convex objectives.
	Some of the modern scalable techniques for the Newton method
	are \textit{block-coordinate updates}, \textit{stochastic subspaces} and \textit{sketching}
	\cite{hanzely2020stochastic,nutini2022let,hanzely2023sketch},
	\textit{distributed} and
	\textit{lazy} computations \cite{qian2021basis,safaryan2021fednl,islamov2021distributed,islamov2022distributed,doikov2023second,doikov2023first}
	and
	\textit{stochastic preconditioning} \cite{pasechnyuk2022effects,pasechnyuk2023convergence},
	for convex and non-convex problems.
	
	Our analysis is based on the new concept of \textit{graded non-convexity},
	which is related to studies of other generalized notions of convexity~\cite{vial1983strong,hormander2007notions}.
	At the same time, more refined specifications of the problem class
	that investigate the distribution of the eigenvalues were considered recently in 
	\cite{kovalev2018stochastic,scieur2020universal,cunha2022only,goujaud2022super}.
	The motivation of our work to cut large gaps between the leading eigenvalues
	is closely related to recently proposed coordinate methods
	with \textit{volume sampling} \cite{rodomanov2020randomized}
	and \textit{polynomial preconditioning} technique \cite{doikov2023polynomial},
	analysing the convex problems with a specific structure.

	\vspace*{-0.5em}
	\paragraph{Contents.}
	The rest of the paper is organized as follows.
	In Section~\ref{SectionProblems} we introduce
	the notion of \textit{graded non-convexity}.
	We study its main properties and provide several examples.
	Section~\ref{SectionSpectral} contains our main algorithm \eqref{MainAlgorithm}
	the Gradient Method with Spectral Preconditioning. We prove 
	fast convergence rates for this method, showing an improvement 
	when increasing the preconditioning order (Theorem~\ref{TheoremLip}).
	In Section~\ref{SectionCut} we show a simple modification of our method
	that allows to remove 
	the dependency on the negative part of the spectrum (Theorem~\ref{TheoremCut}). In Section~\ref{SectionConvex} we show
	improved rates of convergence for the special case
	of convex functions.
	In Section~\ref{SectionImplementation} we discuss the efficient implementation
	of our algorithms.
	Section~\ref{SectionExperiments} presents numerical experiments.
	Missing proofs are provided in the appendix.

	\vspace*{-0.5em}
	\paragraph{Notation.}
    
    We are interested in solving the following unconstrained optimization problem\footnote{Our results can be also generalized to the \textit{composite} formulation of optimization problems, see Section~\ref{SectionAppendixComposite} in the appendix.}
    \beq \label{MainProblem}
    \ba{rcl}
    \min\limits_{\xx \in \R^n} 
    f(\xx) 
    \ea
    \eeq
    where $f: \R^n \to \R$ is a several times differentiable target function, which can be \textit{non-convex}.
    We denote its global lower bound by $f^{\star} := \inf_{\xx} f(\xx)$, which we assume to be finite.
    In non-convex optimization, our goal is to find an approximate \textit{stationary point} $\bar{\xx}$ to \eqref{MainProblem},
    ensuring $\| \nabla f(\bar{\xx}) \| \leq \varepsilon$, for some $\varepsilon > 0$.
    We use the standard Euclidean norm for vectors: $\| \xx \| := \la \xx, \xx \ra^{1/2}, \xx \in \R^n$,
    and the corresponding operator norm for matrices:
    $
    \| \mat{A} \| :=  \max_{\xx : \|\xx\| \leq 1} \| \mat{A} \xx \|.
    $
    
    We denote the gradient of $f$ at point $\xx \in \R^n$ by $\nabla f(\xx) \in \R^n$,
    and the Hessian matrix by $\nabla^2 f(\xx) \in \R^{n \times n}$.
    Note that the Hessian is a symmetric matrix. Hence, for any point $\xx \in \R^n$
    we can introduce the following spectral decomposition:
    \beq \label{SpecDecompose}
    \ba{rcl}
    \nabla^2 f(\xx) & = & \sum\limits_{i = 1}^n \lambda_i(\xx) \uu_{i}(\xx) \uu_{i}(\xx)^{\top},
    \ea
    \eeq
    where the eigenvalues are sorted in a \textit{non-ascending} order:
    $
    \lambda_1(\xx) \geq \lambda_2(\xx) 
    \geq \ldots \geq
    \lambda_n(\xx),
    $
    and $\uu_1(\xx), \ldots, \uu_n(\xx) \in \R^n$ are the corresponding
    orthonormal eigenvectors.
    There are always several possible choices for the eigenvector basis, 
    hence, we assume that a specific selection has been fixed. 
    Our results remain independent of the particular selection. 
    For a fixed spectral decomposition~\eqref{SpecDecompose},
    we denote 
    by $\nabla^2_{\tau} f(\xx)$, $1 \leq \tau \leq n$,
    the \textit{Hessian of spectral order $\tau$}:
    	\vspace*{-0.5em}
    \beq \label{HessOrder}
    \ba{rcl}
    \nabla^2_{\tau} f(\xx)
    & := &
    \sum\limits_{i = 1}^\tau \lambda_i(\xx) \uu_{i}(\xx) \uu_{i}(\xx)^{\top}
    \in \R^{n \times n}.
    \ea
    \eeq
    For convenience, we set $\nabla^2_0 f(\xx) \equiv \0$.
    Thus, $\rank(\nabla^2_{\tau} f(\xx)) \leq \tau$,
    and for  $\tau = n$ we obtain the full Hessian.
    Of course, the decomposition of the form~\eqref{HessOrder} is not unique,
    especially if $\lambda_{\tau}(\xx) = \lambda_{\tau+1}(\xx)$ for a certain $\xx$.
    However, we always assume that a unique choice has been fixed for ease of notation.
    We denote by $\nabla^3 f(\xx)$ the third derivative of $f$ at point $\xx$,
    which is a trilinear symmetric form. The action of this form onto
    some fixed directions $\uu_1, \uu_2, \uu_3 \in \R^n$ is denoted by
    $
    \nabla^3 f(\xx)[ \uu_1, \uu_2, \uu_3 ] \in \R.
    $

	\section{Problem Classes}
	\label{SectionProblems}
	
	A standard assumption in non-convex optimization is to assume that the objective function is smooth, i.e.\ that its gradients are Lipschitz continuous.
    For twice-differentiable functions, this is equivalent to assuming that the norm of the Hessian, $\|  \nabla^2 f(\xx) \|, $ is uniformly bounded.
    In this section, we introduce a new concept that allows us to capture the distribution of the eigenvalues in more detail.
	
	\subsection{Grade of Non-convexity}

	We start with a formal definition of our problem class.
	
	\BD \label{DefinitionGraded}
	For a twice continuously differentiable function $f: \R^n \to \R$
	and convex set $Q \subseteq \R^n$,
	we say that $f$ 
	is \textit{non-convex of grade $\tau$}, if
	\beq \label{GradeNonConvex}
	\boxed{
		\ba{rcl}
		\nabla^2_{\tau} f(\xx) & \succeq & \0, \qquad \forall \xx \in Q.
		\ea
	}
	\eeq
	\ED
	
	In other words, \eqref{GradeNonConvex} means that the top $\tau$ eigenvalues of the Hessian
	are non-negative everywhere on $Q$:
	\beq \label{NonNegativeLambda}
	\ba{rcl}
	\lambda_{\tau}(\xx) & \geq & 0, \qquad \forall \xx \in Q.
	\ea
	\eeq
	In our analysis, we will mostly use for simplicity $Q \equiv \R^n$ 
	(the whole space). However, it can only refine our problem class
	if some localization of a solution $\xx^{\star} \in Q$ is available.
	Definition~\ref{DefinitionGraded} implies a certain 
	restriction on a surface structure of the objective function.
	In differential geometry, 
	condition~\eqref{NonNegativeLambda} on the curvatures leads to the notion
	of \textit{$\tau$-convex surface}~\cite{gromov1991sign}.

	For  $0 \leq \tau \leq n$, we denote by $\mathcal{F}_\tau$ the set of functions $f \in \mathcal{F}_{\tau}$
	that satisfies~\eqref{GradeNonConvex}.
	By our definition, $\F_0 = C^2(\R^n)$ is the set of \textit{all} twice continuously differentiable 
	functions\footnote{
		Note that definition of our classes $\F_{\tau}, 0 \leq \tau \leq n$ depends on set $Q$ and on the problem dimension $n$.
		Thus, it would be more formal to use notation $\F_{\tau}^n(Q)$.
		However, we omit extra indices since they should always be clear from the context.} 
	and $\F_n$ consists only of convex objectives.
	These classes are closed under multiplications by
	a non-negative scalar: $f \in \F_{\tau} \Rightarrow \alpha f \in \F_{\tau}$ for any $\alpha \geq 0$.
	Hence, we obtain the \textit{nested family of functional cones}:
	$$
	\boxed{
		\ba{rcl}
		\underset{\textbf{all functions}}{\F_{0}}  \supset \;\; \F_{1}
		\;\; \supset & \ldots &\supset \;\; \F_{n - 1}
		\;\; \supset 
		\underset{\textbf{convex functions}}{\F_n}
		\ea
	}
	$$
	Intuitively, functions with larger grades should be easier to minimize.
	At the same time, a method that works for a certain $\F_{\tau}$
	can also tackle all problems from $\F_i$ for $i \geq \tau$.

	\subsection{Main Properties}
	\label{SubsectionProperties}
	
	Let us study some of the most basic properties
	that follow from our definition.
	First, we have the following important \textit{grading rule},
	that equips our sets with additional structure.
	
	\BP \label{ProposSum}
	Let $f \in \F_i$ and $g \in \F_j$, for some $0 \leq i, j \leq n$
	such that $i + j \geq n$. Then, it holds:
	$$
	\ba{rcl}
	f + g & \in & \F_{i + j - n}, \\[3pt]
	\smax(f, g) & \in & 
	\F_{i + j - n},
	\ea
	$$
	where
	$\smax(f, g)(\xx) \Def \ln( e^{f(\xx)} + e^{g(\xx)} ) $ is the soft maximum of two functions.
	\EP
	
	In particular, the summation of a function $f \in \F_\tau$ with \textit{any convex} function
	$g \in \F_n$
	cannot decrease its grade~$\tau$. 
	The same holds for (soft) maximum operations\footnote{It holds 
		$\max(f, g) = \lim_{\mu \to 0} \mu \cdot \smax(f / \mu, g / \mu) $.
		We prefer to work with soft max operation to keep all functions in the smooth class.}.
	Next, we observe that the grade is preserved under affine substitutions.
	
	\BP \label{ProposAffine}
	Let $f \in \F_{\tau}(\R^n)$ be non-convex of grade~$\tau$, and let  
	$\mat{A} \in \R^{n \times m}, \bb \in \R^n$, with $m + \tau \geq n$.
	Denote $g(\xx) = f(\mat{A}\xx + \bb)$.
	Then, $g \in \F_{m - n + \tau}(\R^m)$ is non-convex of grade $m - n + \tau$.
	\EP
	
	For $\tau \geq 1$, our Definition~\ref{DefinitionGraded} 
	implies that the Hessian is \underline{not} negative definite
	at any point: $\nabla^2 f(\xx) \not\prec \0$. This condition means that the function 
	\textit{cannot have strict local maxima}.
	This fact can be formalized as follows.
	
	\BP \label{ProposMax}
	Let $f \in \F_{\tau}(Q)$ for $\tau \geq 1$.
	Then, the weak maximum principle holds:
	for any compact $K \subset Q$, the maximum
	is always
	achieved at the boundary,
	$$
	\ba{rcl}
	\max\limits_{\xx \in K} f(\xx) & = &
	\max\limits_{\xx \in \partial K} f(\xx).
	\ea
	$$
	\EP
	
	Finally, let us mention an intuitive geometric description
	of the surface of function $f$ of a certain grade $\tau$.
	
		\begin{figure}[h!]
		\begin{center}
			\includegraphics[scale=1.6]{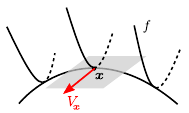}
		\end{center}
		\caption{ \small A surface of a non-convex function $f$.
			At point $\xx$ there is a subspace $V_{\xx}$ where $f$ is convex. }
		\label{FigureSurface}
	\end{figure}

	\BP \label{ProposSurface}
	Let for any $\xx$
	there exists a vector subspace $V_{\xx} \subseteq \R^n$
	with $\dim(V_{\xx}) \geq \tau$ such that
	\beq \label{GeomIneq}
	\ba{rcl}
	f(\xx + \hh) & \geq & f(\xx) + \la \nabla f(\xx), \hh \ra,
	\;\;\; \forall \hh \in V_{\xx}.
	\ea
	\eeq
	Then $f \in \F_{\tau}$ is non-convex of grade $\tau$.
	\EP 
	
	This inequality is stronger than our Definition~\ref{DefinitionGraded}:
	\eqref{GeomIneq} is a sufficient condition for~\eqref{GradeNonConvex}. Geometrically it means
	that for every point $\xx$ there exists a subspace $V_{\xx}$ of the tangent space
	to the surface such that restriction of $f$ onto this subspace is convex (see Fig.~\ref{FigureSurface}).

	\subsection{Examples}
	\label{SubsectionExamples}

	In this section, we provide several examples of non-convex objective functions
	of a non-trivial grade of non-convexity.
	
	\BE[Quadratic Functions] \label{ExampleQuadratic}
	Let $f(\xx) = \frac{1}{2} \la \mat{A} \xx, \xx \ra - \la \bb, \xx \ra$
	for some $\mat{A} = \mat{A}^{\top} \in \R^{n \times n}$ and $\bb \in \R^n$.
	Let the top $\tau$ eigenvalues of $\mat{A}$ be positive:
	$
	\lambda_1(\mat{A}) \geq \ldots \geq  \lambda_{\tau}(\mat{A}) \geq  0
	$.
	Then $f \in \F_{\tau}$.
	\EE
	\vspace{-0.5em}
	Consequently, adding 
	a power of the Euclidean norm as a simple regularizer
	to a non-convex quadratic function as in Example~\ref{ExampleQuadratic}, we obtain
	for $p > 2$ and some $\sigma > 0$:
	\beq \label{RegQuadratic}
	\ba{rcl}
	f(\xx) & = & \frac{1}{2}\la \mat{A} \xx, \xx \ra - \la \bb, \xx \ra + \frac{\sigma}{p}\| \xx \|^p,
	\ea
	\eeq
	which describes a simple family of non-convex objectives that can realize all possible grades $f \in \F_{\tau}$,
	$0 \leq \tau \leq n$,
	while a global solution $\xx^{\star} = \argmin_{\xx} f(\xx)$
	always exists, since $f$ has bounded sublevel sets. 
	The problems of the form~\eqref{RegQuadratic}
	are important in applications to regularized second-order and high-order methods
	\cite{grapiglia2017regularized,nesterov2019implementable,nesterov2022quartic}.

	\BE[Low-rank Vector Fields] \label{ExampleField}
	Let $\varphi: \R \to \R$ be a univariate (possibly non-convex) function,
	and $\uu: \R^n \to \R^n$ be a general differentiable mapping.
	Set
	\beq \label{VectorField}
	\ba{rcl}
	f(\xx) & = & \varphi( \la \uu(\xx), \xx \ra), \qquad \xx \in \R^n.
	\ea
	\eeq
	Assume that there exists $r \leq n - 1$ such that
	for any $\xx$, it holds
	$r \geq  \rank\bigl(  
	\nabla \uu(\xx) + \nabla^2 \uu(\xx)\xx \bigr)$.
	Then, $f: \R^n \to \R$ is non-convex of grade $n - r - 1$.
	\EE
	
\begin{figure}[h!]
		\begin{center}
			\includegraphics[scale=0.40]{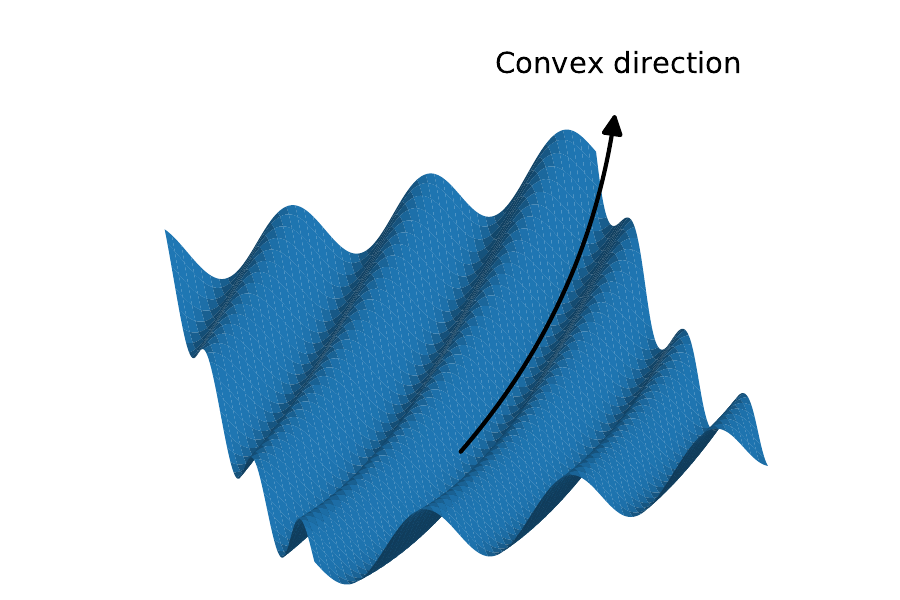}
		\end{center}
		\caption{ \small The graph of two-dimensional function $f(x, y) = \sin ( x + y ) + q(x, y)$, where $q$ is a convex quadratic.
			The non-convex component has the structure of \eqref{VectorField} with $\uu \equiv (1, 1)^{\top}$.  }
		\label{FigSin}
	\end{figure}

	A particular case is when $\uu(\xx) \equiv \uu \in \R^n$ is a constant vector field.
	Then, we obtain that the function \eqref{VectorField}
	is non-convex of grade $n - 1$ (see Fig.~\ref{FigSin}).

	\BE[Partial Convexity] \label{ExamplePartial}
	Let $f(\xx, \yy)$ depend on two groups of variables $\xx \in \R^n$ and $\yy \in \R^m$. Assume that
	for any fixed $\yy$, the function of the first argument
	$$
	\ba{c}
	f(\cdot, \yy) : \; \R^n \to \R 
	\ea
	$$
	is convex. Then $f: \R^{n + m} \to \R$ is non-convex of grade $n$.
	\EE

	We see that if the function is convex with respect to \textit{some} of the variables,
	then, as a function of all variables, it is non-convex of the corresponding grade.
	Note, however, that the actual structure of the subspace when the Hessian is positive can
	be quite complicated.
	As a direct consequence, we obtain that matrix factorizations
	and deep neural network models with convex losses satisfy our assumption.
	
	\BE[Diagonal Neural Networks]  \label{ExampleDNN}
	For a given $\cc \in \R^n$, consider
	\beq \label{DNNEx}
	\ba{rcl}
	f(\xx, \yy) & = & \frac{1}{2} \| \xx \circ \yy - \cc \|^2, \quad \xx, \yy \in \R^n,
	\ea
	\eeq
	where $\circ$ is the element-wise product. At every point $(\xx, \yy) \in \R^{2n}$
	the Hessian has the following set of $2n$ eigenvalues, a pair for each $1 \leq i \leq n$:
	$
	\lambda_i  =  \frac{1}{2}\bigl[ \,  x_i^2 + y_i^2 \pm 
	\sqrt{  (x_i^2 - y_i^2)^2 + 4(2x_i y_i - c_i)^2  } \, \bigr],
	$
	with at least one being always non-negative.
	Therefore, $f: \R^{2n} \to \R$ is non-convex of grade $n$.
	\EE
	 
	The objective of the form \eqref{DNNEx}
	is a good model for studying the dynamics of gradient methods
	in deep learning
	\cite{Woodworth2020,Pesme2021,even2023s,pesme2023saddle}.

	\BE[Matrix Factorizations] \label{ExampleMF}
	For a given target matrix $\mat{C} \in \R^{n \times m}$, 
	and some $r > 0$, consider
	\beq \label{MFEx}
	\ba{rcl}
	f(\mat{X}, \mat{Y}) & = & \frac{1}{2}\| \mat{X} \mat{Y} - \mat{C}  \|^2_F,
	\ea
	\eeq
	where  $\mat{X} \in \R^{n \times r}, \; \mat{Y} \in \R^{r \times m}$.
	Then, $f: \R^{n \times r + r \times m} \to \R$ is non-convex of grade $\max\{ m, n \} \times r$.
	More generally, the function
	$
	f(\mat{X}_1, \ldots, \mat{X}_{d}) =  \frac{1}{2} \| \mat{X}_1 \mat{X}_2 \cdots \mat{X}_d - \mat{C} \|_F^2,
	$
	$\mat{X}_i \in \R^{n_i \times m_i}$,
	is non-convex of grade $\max_{1 \leq i \leq d} \bigl[ n_i \times m_i \bigr]$.
	\EE

	We observe that in many practical scenarios, it is very common to encounter objectives that have 
	a subspace with positive eigenvalues of the Hessian.
	Indeed, the opposite seems rather quite rate and indicates that 
	the target objective is \textit{purely concave}.
	On the contrary, for deep neural network models with convex losses
	we always ensure the existence of a subspace with \textit{positive curvature},
	which serves as the main computational burden for a method
	to converge to a stationary point, when the problem is ill-conditioned. 
	
	In the next sections,
	we propose the spectral preconditioning technique
	for gradient methods in order to tackle ill-conditioning 
	of the positive part of the spectrum. At the same time, as we will show,
	the impact of the negative curvature to optimization methods can easily be alleviated.

	\section{Spectral Preconditioning}
	\label{SectionSpectral}
	
    In this section we present our proposed algorithm.
    The method aims to exploit the (possibly) convex-like structure of the function. 
    When there is no such structure ($\tau = 0$)  the method becomes equal to the standard gradient method.

	At every iteration of our method, we choose 
	a matrix 
	$\mat{H} = \mat{H}^{\top} \succeq 0$
	and perform the following preconditioned gradient step,
	for a given point $\xx \in \R^n$
	and gradient vector $\gg \in \R^n$:
	$$
	\ba{rcl}
	\text{GradStep}_{\mat{H}, \alpha}(\xx, \gg) 
	& := & 
	\xx - \bigl(  \mat{H} + \alpha \mat{I} \bigr)^{-1} \mat{g},
	\ea
	$$
	where $\alpha \geq 0$ is some regularization parameter.
	Hence, the matrix $\mat{H}$ plays the role of a \textit{preconditioner}.
	We want to choose it as an approximation of the Hessian
	of a certain spectral order:
	$\mat{H} \approx \nabla^2_{\tau} f(\xx)$.
	When $\tau = 0$, we have $\nabla^2_{\tau} f(\xx) \equiv \0$
	and thus we do a step of plain gradient descent.
    In contrast, for $\tau = n$ we approximate the full Newton step.
	Let us present the method in algorithmic form.

	\beq\label{MainAlgorithm}
	\ba{|c|}
	\hline\\[-7pt]
	\quad \mbox{\bf Gradient Method with} \quad\\
	\mbox{ \bf{Spectral Preconditioning} } \\
	\\[-7pt]
	\hline\\
	\ba{l}
	\mbox{{\bf Choose} $\xx_0 \in \R^n$ and $0 \leq \tau \leq n$.} \\[5pt]
	\mbox{\bf For $k \geq 0$ iterate:}\\[5pt]
	\mbox{1. Estimate $\mat{H}_k \; \approx \;\nabla^2_{\tau} f(\xx_k) \; \in \; \R^{n \times n} $} \\[5pt]
	\mbox{2. Perform the gradient step, for some $\alpha_k \geq 0$:} \\[5pt]
	$$
	\ba{rcl}
	\;\;
	\xx_{k + 1} & = & 
	\text{GradStep}_{\mat{H}_k, \alpha_k} ( \xx_k, \nabla f(\xx_k) )
	\ea
	$$ \\[5pt]
	\ea\\
	\hline
	\ea
	\eeq
	Up to now, we do not specify explicitly how we estimate $\nabla^2_{\tau} f(\xx_k)$.
	Our matrix $\mat{H}_k$ should be easily computable and we
	aim to maintain a low rank representation,
	\vspace*{-1em}
	\beq \label{HLowRank}
	\ba{rcl}
	\mat{H}_k & = & \sum\limits_{i = 1}^{\tau} a_{k, i} \vv_{k,i} \vv_{k, i}^{\top},
	\ea
	\eeq
	for a set of positive numbers $( a_{k, i} )_{i = 1}^{\tau}$
	and orthonormal vectors $( \vv_{k, i} )_{i = 1}^{\tau}$,
	so that we can perform iterations of algorithm~\eqref{MainAlgorithm} cheaply.
	We discuss details on efficient implementation of every step in Section~\ref{SectionImplementation}.
	To quantify the approximation errors, we denote
	\beq \label{DeltaKDef}
	\ba{rcl}
	\delta_k & := & \| \mat{H}_k - \nabla^2_{\tau} f(\xx_k) \|,
	\quad 
	\delta \; := \; \max\limits_{k} \delta_k. 
	\ea
	\eeq
	Thus, if $\delta = 0$, we use the exact Hessian of spectral order $\tau$.

	\section{Global Convergence}
	\label{SubsectionGlobal}
	
	Let us show the main convergence result for algorithm~\eqref{MainAlgorithm}.
	We establish fast convergence rates to a stationary point
	of our objective~\eqref{MainProblem}, starting from an arbitrary initial point $\xx_0$.
	These rates become better when increasing the spectral order parameter $\tau$ used in our method.
	We introduce the following characteristics of smoothness for our problem classes. 
	We assume that the Hessian is Lipschitz continuous, with parameter $L \geq 0$,
	for all $\xx, \yy$:
	\beq \label{LipHessian}
	\ba{rcl}
	\| \nabla^2 f(\xx) - \nabla^2 f(\yy) \| & \leq & L \|\xx - \yy\|.
	\ea
	\eeq
	To estimate the error of using the Hessian of a spectral order~$\tau$, we
	also introduce the following system of parameters:
	\beq \label{SigmaHessian}
	\ba{rcl}
	\sigma_{\tau} 
	& := &
	\sup_{\xx \in \R^n}
	\| \nabla^2 f(\xx) - \nabla^2_{\tau}  f(\xx) \|.
	\ea
	\eeq
	Thus, $\sigma_0 := \sup_{\xx} \| \nabla^2 f(\xx) \|$ is simply the (best) Lipschitz constant of the gradient of $f$, and $\sigma_n = 0$.
	In general,
	the value of $\sigma_{\tau}$ is the uniform bound for both the $(\tau + 1)th$ eigenvalue of the Hessian,
	and its negative part:
	\beq \label{SigmaBound}
	\ba{rcl}
	\sigma_{\tau} & \geq & 
	\max\{ \lambda_{\tau + 1}(\xx), \, -\lambda_{n}(\xx) \}, 
	\quad \forall \xx \in \R^n.
	\ea
	\eeq
	Under these conditions, and using a sufficiently big regularization parameter $\alpha_k$ 
	(a second-order ``step size"),
	we can prove the following progress for one iteration $\xx_k \mapsto \xx_{k + 1}$:
	\BL \label{LemmaProgress}
	Let $\alpha_k \geq \sqrt{\frac{L \| \nabla f(\xx_k) \|}{2}} + \sigma_{\tau} + \delta$. Then
	\beq \label{MainFuncProgress2}
	\ba{rcl}
	f(\xx_k) - f(\xx_{k + 1}) & \geq &
	\frac{1}{8 \alpha_k} \| \nabla f(\xx_{k + 1})\|^2. 
	\ea
	\eeq
	\EL
	\vspace*{-1em}
	Note that we do not need to know the exact values of the parameters $L$, $\sigma_{\tau}$, and $\delta_k$. In practice, we can use an adaptive search 
	to ensure sufficient progress \eqref{MainFuncProgress2}.
	This lemma leads us to the basic global convergence result.
	\BTF \label{TheoremLip}
	Let $f \in \F_{\tau}$ be non-convex of grade~$\tau$,
	where $0 \leq \tau \leq n$ is fixed.
	Let $f$ have a Lipschitz Hessian with constant $L$
	and bounded parameter $\sigma_{\tau}$. Consider iterations $\{ \xx_k \}_{k \geq 0}$
	of algorithm~\eqref{MainAlgorithm} with
	\vspace{-0.5em}
	\beq \label{Alg1AlphaK}
	\ba{rcl}
	\alpha_k & = & \sqrt{\frac{L \| \nabla f(\xx_k) \|}{2}} + \sigma_{\tau} + \delta.
	\ea
	\eeq
	Then, for any $\varepsilon > 0$, it is enough to do $K = $
	\vspace{-0.5em}
	$$
	\ba{cl}
	&
	\!\!\!\!\!\!
	\!\!\!\!\!\!
	\bigl\lceil
	8 ( f(\xx_0) - f^{\star} ) \cdot \bigl(   \sqrt{\frac{L}{2}} \frac{1}{\varepsilon^{3/2}} 
	+ \frac{\sigma_{\tau} + \delta}{\varepsilon^2}  \bigr) 
	\; + \; 2 \ln \frac{\| \nabla f(\xx_0) \|}{\varepsilon}
	\bigr\rceil
	\ea
	$$
	steps to ensure $\min\limits_{1 \leq i \leq K} \| \nabla f(\xx_i) \| \leq \varepsilon$.
	\ETF

	For $\tau = n$ (the full Newton), we have $\sigma_n = 0$
	and the complexity becomes $\cO( \frac{1}{\varepsilon^{3/2}} ) $
	up to an additive logarithmic term.
	It can also be established for the cubic regularization of the Newton method \cite{nesterov2006cubic}.
	For $\tau = 0$ we recover the rate of the gradient descent on non-convex problems \cite{nesterov2018lectures}.
	Note that 
	the rule \eqref{Alg1AlphaK} for choosing 
	$\alpha_k$  is very simple.
	It is inspired by the gradient regularization technique
	developed initially for convex optimization
	\cite{mishchenko2023regularized,doikov2023gradient}.

	\section{Cutting the Negative Spectrum}
	\label{SectionCut}

	In our previous result, we saw that the complexity of the method
	depends on the parameter $\sigma_{\tau}$, which becomes
	smaller when increasing the spectral order $\tau$
	of the method.
	However, due to \eqref{SigmaBound}, it includes a bound
	for the absolute value of the negative curvature, which can be constantly big.
	In this section, we propose a simple modification of our step-size rule,
	which alleviates this issue. 
	
	Let us denote the \textit{positive part} of the Hessian, 
	$\nabla^2_+ f(\xx) :=  \nabla^2_{\xi} f(\xx)$,
	where
	$\xi = \argmax\{  \tau : \nabla^2_{\tau} f(\xx) \succeq \mat{0} \}$.
	Correspondingly, the \textit{negative part}
	is $\nabla^2_{-}f(\xx) := \nabla^2_{+}f(\xx) - \nabla^2 f(\xx) \succeq \mat{0}$.
	We introduce the following system of parameters (compare with \eqref{SigmaHessian}):
	\beq \label{SigmaPlus}
	\ba{rcl}
	\sigma_{\tau}^+ 
	& := & 
	\sup\limits_{\xx \in \R^n}
	\| \nabla^2_{+} f(\xx) - \nabla^2_{\tau} f(\xx) \|.
	\ea
	\eeq
	So, $\sigma_{\tau}^+$ is the uniform bound
	for the positive tail of the spectrum,
	and it is no longer affected by the negative curvature:
	\beq \label{SigmaPlusBound}
	\ba{rcl}
	\sigma_{\tau}^+ & \geq & 
	\max\{ \lambda_{\tau + 1}(\xx), \, 0 \}, 
	\quad \forall \xx \in \R^n.
	\ea
	\eeq
	Our new step-size rule is based on the cubic regularization technique. 
	Namely, at every iteration $k \geq 0$, we compute the regularization 
	parameter as the maximum of the following univariate concave function:
    $$
	\ba{rcl}
	\alpha_k^{\star} 
	& \!\!\!\! = \!\!\!\! & 
	\argmax\limits_{\alpha > 0}\Bigl[ -\frac{1}{2} \la (\mat{H}_k + (\alpha + \eta) \mat{I})^{-1} \gg_k, \gg_k \ra
	- \frac{2 \alpha^3}{3L^2}\Bigr],
	\ea
	$$
	where 
	$\gg_k := \nabla f(\xx_k)$ is the current gradient, 
	$\eta \geq 0$ is a balancing term,
	and $L > 0$ is the Lipschitz constant of the Hessian.
	This subproblem is well defined
	and has a unique global maximum, which can be found by using
	binary search or univariate Newton iterations \cite{conn2000trust,nesterov2006cubic}.
	Then, to eliminate the negative curvature, 
	for the sequence generated by our method
	we use an extra sequence of test points $\{ \yy_k \}_{k \geq 1}$
	defined by
	$\yy_{k + 1}  := \xx_k + \mat{P}_k(\xx_{k + 1} - \xx_k)$,
	where $\mat{P}_k$ is a projector which preserves the image of $\mat{H}_k$:	
	$\mat{H}_k \mat{P}_k =  \mat{H}_k$,
	but has a small intersection with the negative part of the Hessian,
	$$
	\ba{rcl}
	\!
	\| \nabla^2_{-} f(\xx_k) \mat{P}_k  \| 
	& \!\!\!\! = \!\!\!\! &
	\|  ( \nabla^2 f(\xx_k) - \nabla^2_{+} f(\xx_k) ) \mat{P}_k \|
	\; \leq \; \delta_{-},
	\ea
	$$
	for some $\delta_{-} > 0$.
	It can be built directly from the low-rank representation~\eqref{HLowRank}
	of $\mat{H}_k$, as follows:
	$\mat{P}_k = \sum_{i = 1}^{\tau} \vv_{k, i} \vv_{k, i}^{\top}$.
	Therefore, this matrix comes with no extra cost. We have
	$\delta_{-} = 0$ when $\vv_{k, i}$ are orthogonal
	to all eigenvectors of $\nabla^2_{-} f(\xx_k)$.
	Note that the tolerance parameters $\delta$ and $\delta_{-}$ are induced 
	by the approximation error $\mat{H}_k \approx \nabla^2_{\tau}f(\xx_k)$,
	and, contrary to $\sigma_{\tau}$ and $\sigma_{\tau}^+$, they do not describe
	our problem class.
	Using a high accuracy approximations for matrices $\mat{H}_k$ (see Section~\ref{SectionImplementation}),
	we can make both $\delta$ and $\delta_{-}$ arbitrarily close to zero.
	We are ready to state the better complexity result 
	for algorithm~\eqref{MainAlgorithm}
	with a convergence rate
	that is independent of the negative curvature.
	\BTF \label{TheoremCut}
	Let $f \in \F_{\tau}$ be non-convex of grade~$\tau$,
	where $0 \leq \tau \leq n$ is fixed.
	Let $f$ have a Lipschitz Hessian with constant $L$
	and bounded parameter $\sigma_{\tau}^+$. Consider iterations $\{ \xx_k \}_{k \geq 0}$
	of algorithm~\eqref{MainAlgorithm} with 
	\beq \label{AlphaCubicChoice}
	\ba{rcl}
	\alpha_k & = & \alpha_k^{\star} + \eta,
	\ea
	\eeq
	where $\eta := \sigma_{\tau}^+ + \delta + \delta_{-}$ 
	is the balancing term.
	Then, for any $\varepsilon > 0$, it is enough to do
	$$
	\ba{cl}
	&
	\!\!\!\!\!\!
	\!\!\!\!\!\!
	\!
	K \; = \;
	\Bigl\lceil 
	2 ( f(\xx_0) - f^{\star} ) \Bigl(   \frac{3\sqrt{2L}}{\varepsilon^{3/2}} 
	+ \frac{16\eta}{\varepsilon^2}  \Bigr) \!
	\Bigr\rceil
	\ea
	$$
	steps to ensure $\min\limits_{1 \leq i \leq K} \| \nabla f(\yy_i) \| \leq \varepsilon$.
	\ETF
	This complexity bound is similar to that one of Theorem~\ref{TheoremLip},
	where we substituted $\sigma_{\tau} \mapsto \sigma_{\tau}^+$.
	It can be much better in case of large negative eigenvalues of the Hessian, since 
	then $\sigma_{\tau}^+$ is much smaller than $\sigma_{\tau}$.
	The cost of such an improvement
	is a more complex rule \eqref{AlphaCubicChoice}, involving $\alpha_k^{\star}$,
	and the guarantee is given for the auxiliary points $\{ \yy_k \}_{k \geq 1}$.

	\section{Convex Convergence Analysis} 
	\label{SectionConvex}

	In this section, we provide the analysis of our method
	for a specific case of $f \in \F_n$ (convex objectives).
	In this case, we do not have the negative curvature part: $\nabla^2_{-}f(\xx) \equiv \0$,
	and
	$\sigma_{\tau}^+ \equiv \sigma_{\tau}$.
	Since $\F_n \subset \F_{\tau}$ for any $\tau$,
	the results of the previous sections can be applied directly
	for the convex case. However,
	as we show, using a more refined technique, we can
	prove much better convergence rates for convex and strongly convex problems.
	We denote by $\mu \geq 0$ the parameter of strong convexity,
	such that 
	$ \nabla^2 f(\xx) \succeq \mu \mat{I}, \forall \xx$.
	
	When the Hessian is positive semidefinite, we can naturally use it define
	the \textit{local norm} at point $\xx$ by $\| \uu \|_{\xx} := \la \nabla^2 f(\xx) \uu, \uu \ra^{1/2}$.
	The local norm becomes more appropriate for describing 
	the right second-order geometry of the objective \cite{nesterov1994interior}.
	With its help, it is possible to characterize smoothness more accurately. 
	Now, we assume
	that $f$ is \textit{quasi-self-concordant}\footnote{Note that for the functions with Lipschitz continuous Hessian~\eqref{LipHessian},
	we can bound the third derivative with the product of Euclidean norms: 
	$\nabla^3 f(\xx)[\uu, \uu, \vv] \leq L \| \uu \|^2 \| \vv \|$, 
	which is in most cases less accurate than \eqref{MquasiSC}.} with parameter $M \geq 0$,
	for all $\xx, \uu, \vv \in \R^n$:
	\beq \label{MquasiSC}
	\ba{rcl}
	\nabla^3 f(\xx)[ \uu, \uu, \vv ] & \leq & M \| \uu \|_{\xx}^2 \| \vv \|.
	\ea
	\eeq
	This condition was considered in
	\cite{bach2010self,karimireddy2018global,sun2019generalized,doikov2023minimizing}.
	It appears that many practical problems satisfy this assumption, including
	\textit{Logistic Regression}, 
	\textit{Soft Maximum}, 
	\textit{Matrix Balancing}, and 
	\textit{Matrix Scaling} problems. 
	We adapt algorithm~\eqref{MainAlgorithm} to this problem class, 
	establishing fast convergence rates that improve with increase of the spectral order~$\tau$.
	We denote by $D := \sup \{ \| \xx - \xx^{\star} \| \, : \, f(\xx) \leq f(\xx_0)   \}$ the 
	diameter of the initial sublevel set which we assume to be bounded. We establish the following result. 
	
	\BTF \label{TheoremConvex}
	Let $f \in \F_{n}$ be strongly convex with $\mu > 0$ and quasi-self-concordant with constant $M$.
	Consider iterations $\{ \xx_k \}_{k \geq 0}$ of algorithm~\eqref{MainAlgorithm}
	with
	\beq \label{AlphaQSC}
	\ba{rcl}
	\alpha_k & = &  M \| \nabla f(\xx_k) \| + \sigma_\tau + \delta.
	\ea
	\eeq
	Then, for any $\varepsilon > 0$, it is enough to do $K = $
	$$
	\ba{cl}
	&
	\!\!\!\!\!\!
	\!\!\!\!\!\!
	4 \Bigl\lceil
	\Bigl( MD + \frac{\sigma_{\tau} + \delta}{2\mu}  \Bigr)
	\ln \frac{f(\xx_0) - f^{\star}}{\varepsilon} 
	\; + \; \ln \frac{\| \nabla f(\xx_0)\| D }{\varepsilon}
	\Bigr\rceil
	\ea
	$$
	steps to ensure $f(\xx_K) - f^{\star} \leq \varepsilon$.
	\ETF
	
	\BR
	Note that by adding the quadratic regularizer to our objective 
	with small $\mu := \frac{\varepsilon}{D^2}$,
	we can turn any convex problem into strongly convex one. Hence,
	we obtain the following complexity for the general convex case:
	$$
	\ba{c}
	\cO\Bigl(   \frac{(\sigma_{\tau} + \delta) D^2}{\varepsilon}
	+ MD \ln \frac{f(\xx_0) - f^{\star}}{\varepsilon}  + \ln \frac{\| \nabla f(\xx_0) \| D}{\varepsilon} \Bigr). 
	\ea
	$$
	\ER
	\vspace*{-1em}
	
	Let us consider the simplest example of minimizing a convex quadratic function ($M = 0$),
	using $\tau := 1$.
	Then, according to Theorem~\ref{TheoremConvex}, we need
	$\tilde{\cO}\bigl( \frac{\lambda_2}{\lambda_n} \bigr)$ outer steps.
	Taking into account iterations of the power method, the total complexity becomes
	$\tilde{\cO}\bigl( \frac{\lambda_2}{\lambda_n} \cdot \frac{\lambda_1}{\lambda_1 - \lambda_2} \bigr)$.
	Note that it can be much better than $\tilde{\cO}\bigl( \frac{\lambda_1}{\lambda_n} \bigr)$
	of the basic gradient method, when $\lambda_1 \ggg \lambda_2$.

	\section{Efficient Implementation}
	\label{SectionImplementation}

    The method \eqref{MainAlgorithm} relies on an approximation $\mat{H}_k \approx \nabla^2_{\tau} f(\xx_k)$. 
    We acknowledge that it might be costly to find such an approximation in general. 
    Below, we outline a method that iteratively constructs low rank approximations by employing $T_k \geq 0$ steps of the power method.
    Here $T_k$ is our parameter. 
    In our experiments we use $T_k=1$ (combined with hot-start from $\mat{H}_{k-1}$). However, other options are to choose $T_k=0$ for some of the iterations (no update), or to spend more iteration in the initialization phase. We leave the explorations of such schedules to future work.

	\paragraph{Power Method.}
	
	For our experiments, we use a generalization of the classical power method
	for finding the top $\tau$ eigenvectors of a given matrix,
	which is also known as \textit{orthogonal iterations}. 
	We can write our matrix as
	$\mat{H}_k = \mat{V}_k \Diag(\aa_k) \mat{V}_k^{\top}$,
	where $\aa_k \in \R^{\tau}_{>0}$
	and $\mat{V}_k \in \R^{n \times \tau}$ consists of orthonormal columns $\{ \vv_{k, i} \}_{i = 1}^{\tau}$.
	We denote by $[ \mat{A} ]_{Q} := \mat{Q}$
	the orthogonal matrix from the resulting QR-factorization of $\mat{A} = \mat{Q}\mat{R}$.
	It can be implemented as the standard Gram–Schmidt
	orthogonalization process with arithmetic complexity $\cO(\tau^2 n)$ operations.
	Then, for each 
	$k \geq 0$
	we use the following procedure to update matrix $\mat{V}_k$:
	\beq\label{PowerMethod}
	\ba{|c|}
	\hline\\[-7pt]
	\quad \mbox{\bf Power Method} \quad\\
	\\[-7pt]
	\hline\\
	\ba{l}
 	\mbox{Set $\mat{V}_k := \mat{V}_{k - 1}$ if $k \geq 1$ else random init.} \\[5pt]
	\mbox{\bf For $t = 1 \ldots T_k$ iterate:}\\[7pt]
	\mbox{\qquad Update $\mat{V}_k := \bigl[ \nabla^2 f(\xx_k) \mat{V}_k \bigr]_Q$ } \\[5pt]
	\mbox{Set $a_{k, i} := \la \nabla^2 f(\xx_k) \vv_{k, i}, \vv_{k, i} \ra$.}\\[3pt]
	\ea\\
	\hline
	\ea
	\eeq
	This procedure converges with linear rate, with the main complexity factor
	proportional to the $\tau$-th spectral gap (see, e.g.,Theorem 8.2.2 from \cite{Golub2013}).
	More advanced approaches include Oja's and Lanczos iterations.

	\vspace*{-1em}
	\paragraph{Low-Rank Updates.}
	The Woodbury matrix identity provides us with the following formula:
	\vspace*{-0.2em}
	$$
	\ba{rcl}
	\!\!\!
	( \mat{H}_k + \alpha_k \mat{I} )^{-1} 
	& \!\!\!\!\! = \!\!\!\!\!&
	\frac{1}{\alpha_k}\bigl[  \mat{I} 
	- \!\mat{V}_k \bigl(  \mat{I} + \alpha_k \Diag( \aa_k )^{-1}  \bigr)^{\!-1} \mat{V}_k^{\top} \bigr].
	\ea
	$$
	Therefore, we need only $O( \tau^2 n)$ arithmetical operations to perform the step,
	which is linear with respect to $n$
	and can be implemented very efficiently for small $\tau$.

	\section{Experiments}
	\label{SectionExperiments}

\begin{figure}[h!]
		\begin{center}
			\includegraphics[scale=0.38]{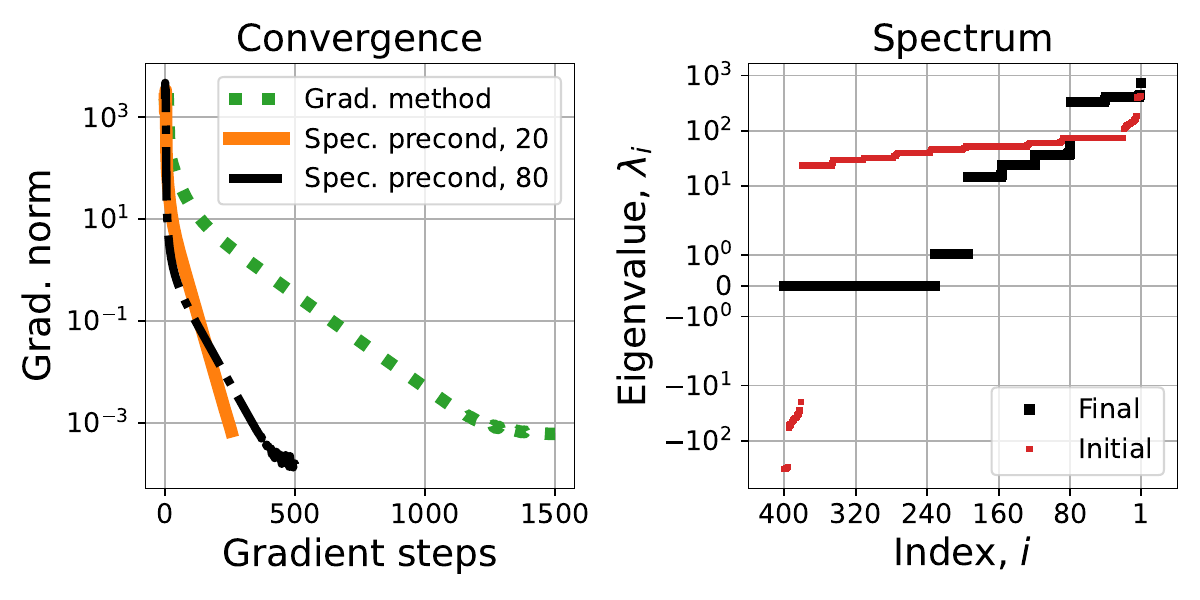}
		\end{center}
		\vspace*{-1em}
		\caption{ \small Matrix factorization, $n = 400$. Left: convergence of the methods with spectral preconditioning for $\tau = 20, 80$. Right: spectrum of the Hessian
			at the initial point and after the last step.
		}
		\label{FigMatFact}
			\vspace*{-1em}
	\end{figure}

	We present illustrative numerical experiments on several
	machine learning problems.
	In Fig.~\ref{FigMatFact}, we study the convergence of the outer iterations of algorithm~\eqref{MainAlgorithm}.
	See Section~\ref{SectionExtraExperiments} in the appendix for the total number of arithmetic operations,
	taking the cost of the power method into account.
		\vspace*{-1em}

	\section{Conclusion}
	\label{SectionConclusion}
	
	In this work, we propose to use inexact top eigenvectors
	of the Hessian as a preconditioning for the gradient methods.
	We introduce a notion of \textit{graded non-convexity},
	which equips our problem classes with a uniform guarantee
	on the number of positive eigenvalues. We show that using
	the preconditioner of order $\tau \geq 1$ provably improves
	the rate of convergence, cutting the gap between
	the top $\tau$ eigenvalues.

	\section*{Acknowledgements}

	This work was supported by the Swiss State Secretariat for
	Education, Research and Innovation (SERI) under contract
	number 22.00133.

	\section*{Impact Statement}
	
	This paper presents work whose goal is to advance the field of Machine Learning. 
	There are many potential societal consequences of our work, none which we feel must be specifically highlighted here.


	\bibliography{bibliography}
	\bibliographystyle{icml2024}

	\newpage
	\appendix
	\onecolumn
	\icmltitle{Appendix}
	
	\section{Extra Experiments}
	\label{SectionExtraExperiments}

\begin{figure}[h!]
		\begin{center}
			\includegraphics[width=0.50\linewidth]{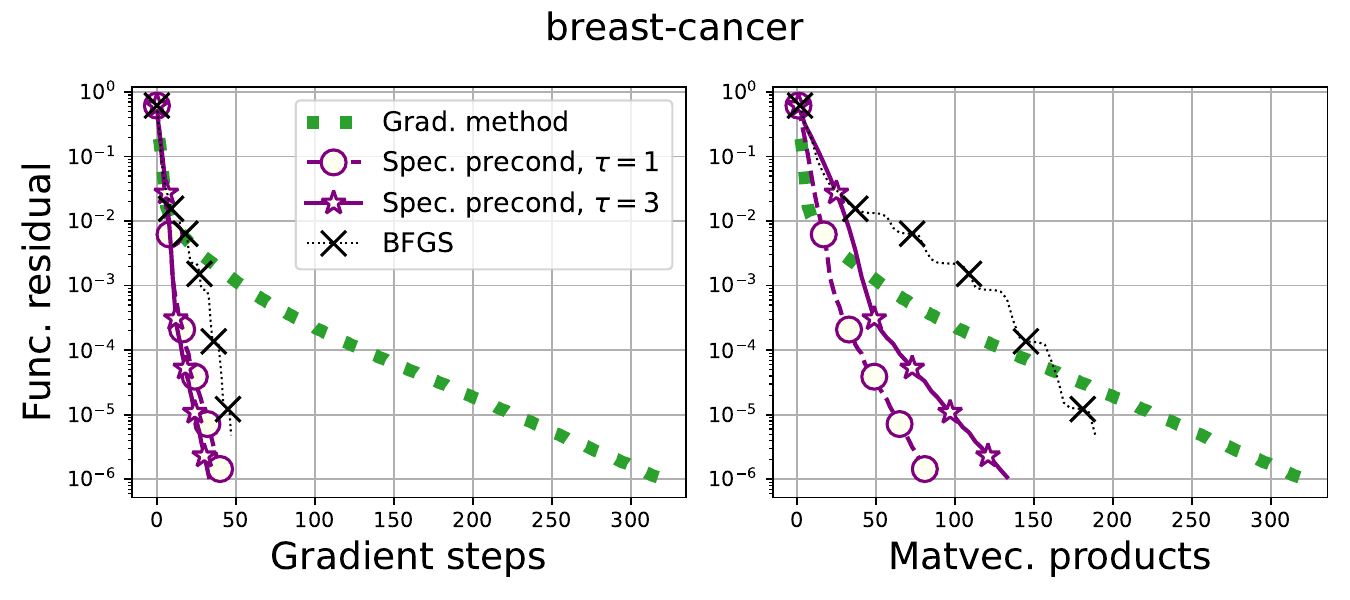}
			\includegraphics[width=0.50\linewidth]{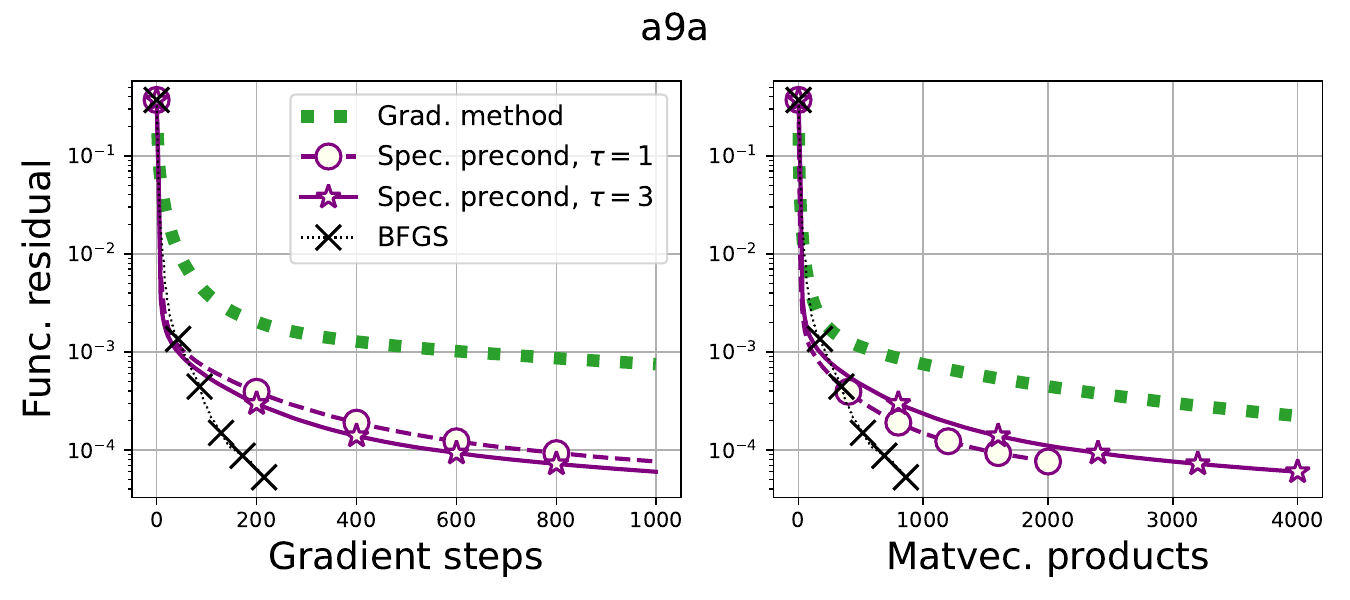}
			\includegraphics[width=0.50\linewidth]{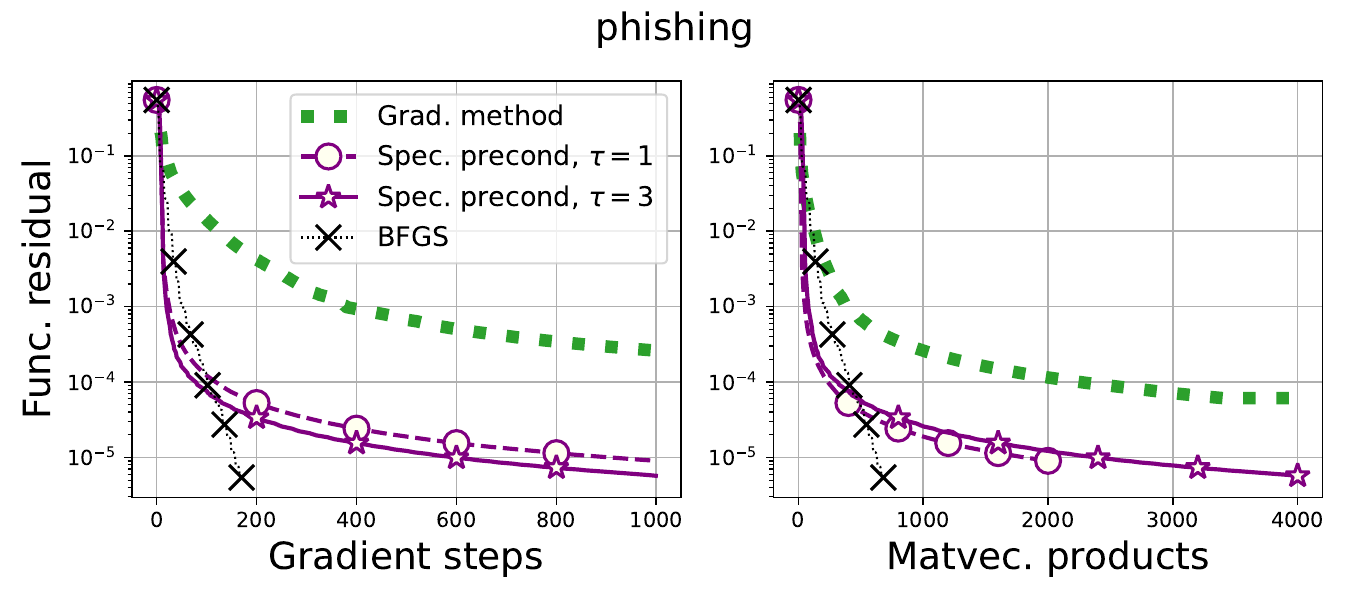}
			\includegraphics[width=0.50\linewidth]{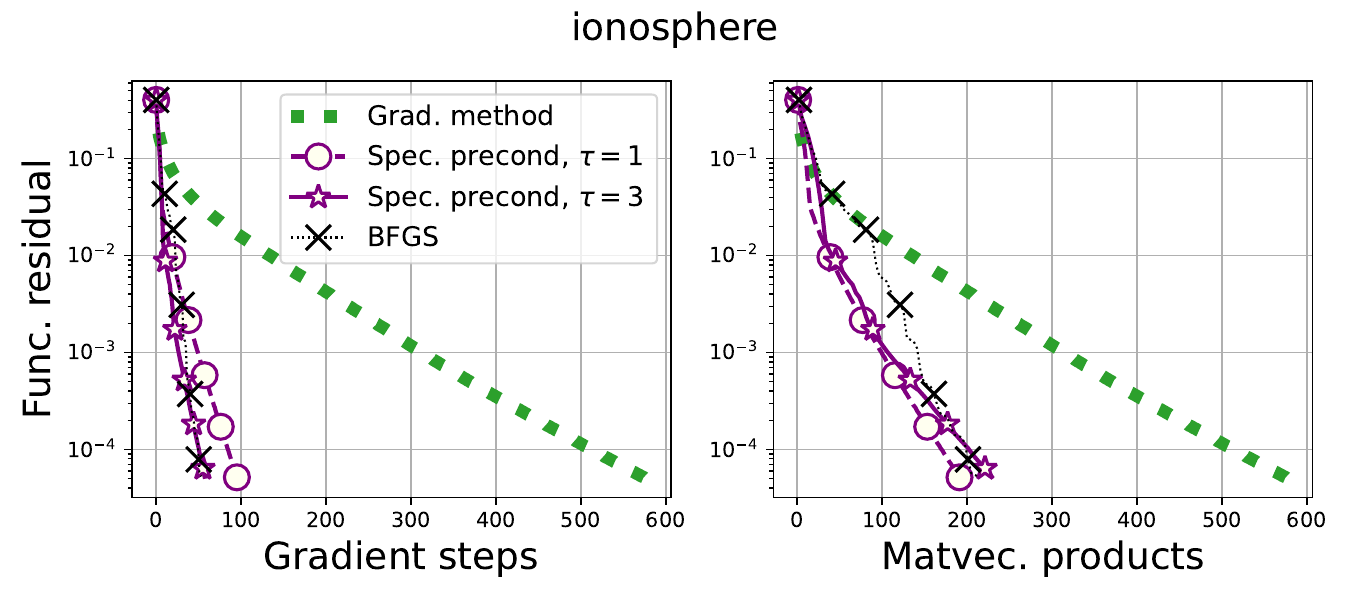}
		\end{center}
		\caption{ \small Training logistic regression with algorithm~\eqref{MainAlgorithm}
			for $\tau = 1$ and $\tau = 3$.
			We see that the spectral preconditioning works significantly better 
			than the basic gradient descent and comparable with
			the powerful BFGS algorithm.
		}
		\label{FigLogReg}
	\end{figure}

	In the following numerical experiments we train a logistic regression
	model on several machine learning datasets\footnote{\url{https://www.csie.ntu.edu.tw/~cjlin/libsvmtools/datasets/}}.
	The results are shown in Fig.~\ref{FigLogReg}.

	\newpage
	\section{Proofs for Section~\ref{SectionProblems}}
	
	\subsection{Proof of Proposition~\ref{ProposSum}}
	
	\proof 
	For any two symmetric matrices $\mat{A}, \mat{B} \in \R^{n \times n}$, 
	and $1 \leq i, j \leq n$ such that $i + j \geq n + 1$, we have
	the following Weyl's inequality:
	\beq \label{EigsSum}
	\ba{rcl}
	\lambda_{i + j - n}(\mat{A} + \mat{B}) & \geq &
	\lambda_i(\mat{A}) + \lambda_j(\mat{B}),
	\ea
	\eeq
	which immediately proves the statement about the sum.
	To prove the statement for soft maximum, 
	let us denote 
	$$
	\ba{rcl}
	f(\xx) & := &
	\smax(f_1(\xx), f_2(\xx))
	\;\; = \;\; \ln( e^{ f_1(\xx) } + e^{ f_2(\xx) } ),
	\ea
	$$
	and compute the derivatives. We obtain, for any $\hh \in \R^n$:
	$$
	\ba{rcl}
	\la \nabla f(\xx), \hh \ra & = & 
	\sigma( f_1(\xx) - f_2(\xx) ) \cdot \la \nabla f_1(\xx), \hh \ra
	+ 
	\sigma( f_2(\xx) - f_1(\xx) ) \cdot \la \nabla f_2(\xx), \hh \ra,
	\ea
	$$
	where $\sigma(t) \Def \frac{1}{1 + e^{-t}}$.
	Note that $\sigma'(t) = \sigma(t) \cdot \sigma(-t)$.
	Denoting $\sigma_1 \Def \sigma( f_1(\xx) - f_2(\xx) )$,
	$\sigma_2 \Def \sigma( f_1(\xx) - f_2(\xx) )$,
	we get
	$$
	\ba{rcl}
	\la \nabla^2 f(\xx) \hh, \hh \ra
	& = & 
	\sigma_1  \sigma_2  \Bigl[  \la \nabla f_1(\xx), \hh \ra - \la \nabla^2 f(\xx), \hh \ra  \Bigr]^2
	+ \sigma_1 \la \nabla^2 f_1(\xx) \hh, \hh \ra
	+ \sigma_2 \la \nabla^2 f_2(\xx) \hh, \hh \ra \\
	\\
	& \geq & 
	\sigma_1 \la \nabla^2 f_1(\xx) \hh, \hh \ra
	+ \sigma_2 \la \nabla^2 f_2(\xx) \hh, \hh \ra.
	\ea
	$$
	Using \eqref{EigsSum} completes the proof.
	\qed

	\subsection{Proof of Proposition~\ref{ProposAffine}}
	
	\proof
	The Hessian of $g$ is given by
	$$
	\ba{rcl}
	\nabla^2 g(\xx) & = & \mat{A}^{\top} \nabla^2 f( L(\xx) ) \mat{A}
	\;\; \in \;\; \R^{m \times m}.
	\ea
	$$
	The linear map $\mat{A}: \R^m \to \R^n$ induces the isomorphism: 
	$\R^m / \ker(\mat{A}) \overset{\varphi}{\to} \im(\mat{A})$ and thus $m  = \dim \ker (\mat{A}) + \dim \im(\mat{A})$ (the classical Rank–nullity theorem).
	
	By our assumption, $\nabla_{\tau}^2 f(L(\xx)) \succeq \0$. Let us denote by $V_+$ the subspace spanned by top $\tau$ eigenvectors of $\nabla^2 f(L(\xx))$, and $V_-$ is spanned by the rest. Hence, $\R^n = V_+ \oplus V_-$. Then
	$$
	\ba{rcl}
	m & = & \dim \ker(\mat{A}) + \dim \im(\mat{A}) \cap V_+
	+ \dim \im(\mat{A}) \cap V_- \\
	\\
	& \leq & 
	\dim \ker(\mat{A}) + \dim \im(\mat{A}) \cap V_+ +  n - \tau.
	\ea
	$$
	Rearranging the terms, we get
	$$
	\ba{rcl}
	\dim \ker(\mat{A}) + \dim \im(\mat{A}) \cap V_+ & \geq & m - n + \tau.
	\ea
	$$
	Therefore, we conclude that the linear subspace $U = \ker(\mat{A}) \, \oplus \, \varphi^{-1}( \im(\mat{A}) \cap V_+ )$
	has dimension $\dim U \geq m - n + \tau$, and for any $\hh \in U$:
	$$
	\ba{rcl}
	\la \nabla^2 g(\xx) \hh, \hh \ra & \geq & 0,
	\ea
	$$
	which proves the required bound.
	\qed
	
	\subsection{Proof of Proposition~\ref{ProposMax}}
	
	\proof 
	For a small $\varepsilon > 0$, we consider the regularized objective
	$$
	\ba{rcl}
	g(\xx) & := & f(\xx) + \frac{\varepsilon}{2}\|\xx\|^2.
	\ea
	$$
	Assume that $g$ achieves its maximum $\xx^{\star}_g = \arg\max\limits_{\xx \in K} g(\xx)$
	in the interior, $\xx^{\star}_g \in \inter K$.
	Then, the second-order stationary condition implies that
	$$
	\ba{rcl}
	\nabla^2 g(\xx^{\star}_g) & \preceq & \mat{0},
	\ea
	$$	
	which is impossible due to $\lambda_1(  \nabla^2 g(\xx^{\star}_g) ) = \lambda_1(\nabla^2 f(\xx^{\star}_g)) + \varepsilon > 0$. 
	Hence, $\xx^{\star}_g \in \partial K$ and we have that
	$$
	\ba{rcl}
	\max\limits_{\xx \in K} {f(\xx)}
	& \leq & 
	\max\limits_{\xx \in K} g(\xx)
	\;\; = \;\;
	\max\limits_{\xx \in \partial K} g(\xx)
	\;\; \leq \;\;
	\max\limits_{\xx \in \partial K} f(\xx) + \frac{\varepsilon}{2} R^2,
	\ea
	$$
	where $R$ is the radius of a ball containing $K$.
	Tending $\varepsilon$ to $0$ completes the proof.
	\qed

	\subsection{Proof of Proposition~\ref{ProposSurface}}
	
	\proof 
	Indeed, by the Taylor theorem, we have that for any $\hh \in V_{\xx}$ and $\lambda > 0$
	there exists $\xi \in [0, 1]$ such that
	$$
	\ba{rcl}
	0 & \leq & f(\xx + \lambda \hh) - f(\xx) - \lambda \la \nabla f(\xx), \hh \ra \\
	\\
	& = &
	\frac{\lambda^2}{2} \la \nabla^2 f(\xx + \lambda \xi \hh) \hh, \hh \ra.
	\ea
	$$	
	Dividing both sides by $\lambda^2$ and
	taking the limit $\lambda \to 0$, we obtain that
	$$
	\ba{rcl}
	\la \nabla^2 f(\xx) \hh, \hh \ra & \geq & 0, \qquad \forall \hh \in V_{\xx}.
	\ea
	$$
	Therefore, since $\dim(V_{\xx}) \geq \tau$, we have that $\lambda_{\tau}(\xx) \geq 0$.
	\qed

	\section{Composite Optimization and Proofs for Section~\ref{SectionSpectral}}
	\label{SectionAppendixComposite}
	
	\subsection{Composite Formulation}
	
	The main results of our work can be generalized to 
	a more broad family of Composite Optimization Problems of the following form:
	\beq \label{MainProblemComposite}
	\ba{rcl}
	\min\limits_{\xx \in Q} 
	\Bigl[ 
	F(\xx) & := & f(\xx) + \psi(\xx)
	\Bigr],
	\ea
	\eeq
	where $f : \R^n \to \R$ is our smooth objective (which can be non-convex),
	and $\psi: \R^n \to \R \cup \{ +\infty \}$ is a \textit{simple} closed convex function
	(e.g., indicator of a given closed convex set, or a regularizer). We set $Q := \dom \psi \subseteq \R^n$.
	Thus, the main properties of $f$ (the grade of non-convexity and the level of smoothness
	can be defined with respect to this convex set $Q$).
	
	In case of presence of the composite part $\psi(\cdot)$ in our problem, iterations of our method should be modified.
	For a given point $\xx \in Q$, gradient vector $\gg \in \R^n$ and matrix $\mat{H} = \mat{H}^{\top} \succeq 0$,
	we define the composite gradient step, as follows:
	\beq \label{CompositeStep}
	\ba{rcl}
	\text{CompositeStep}_{\mat{H}, \alpha}(\xx, \gg)
	& := &
	\argmin\limits_{\yy \in Q}
	\Bigl\{\,
	\la \gg, \yy \ra + \psi(\yy) +  \frac{1}{2} \la (\mat{H} + \alpha \mat{I}) (\yy - \xx), \yy - \xx \ra
	\Bigr\},
	\ea
	\eeq
	where $\alpha \geq 0$ is the regularization parameter.
	In general, due to the presence of regularization term in \eqref{CompositeStep}, for $\alpha > 0$,
	the composite subproblem is \textit{strongly convex}, and we can 
	employ the fast linear convergence of gradient methods as applied to \eqref{CompositeStep},
	for computing this step inexactly.

	With these modifications, we are ready to present a composite version of our algorithm
	for solving \eqref{MainProblemComposite}:
	
	\beq\label{CompositeMethod}
	\ba{|c|}
	\hline\\[-5pt]
	\quad \mbox{\bf Composite Gradient Method with} \quad\\
	\mbox{ \bf{Spectral Preconditioning} } \\
	\\[-5pt]
	\hline\\
	\ba{l}
	\mbox{{\bf Choose} $\xx_0 \in Q$ and $0 \leq \tau \leq n$.} \\[7pt]
	\mbox{\bf For $k \geq 0$ iterate:}\\[7pt]
	\mbox{1. Estimate $\mat{H}_k \; \approx \;\nabla^2_{\tau} f(\xx_k) \; \in \; \R^{n \times n} $} \\[7pt]
	\mbox{2. Perform the composite gradient step, for some $\alpha_k \geq 0$:} \\[7pt]
	$$
	\ba{rcl}
	\;\;
	\xx_{k + 1} & = & 
	\text{CompositeStep}_{\mat{H}_k, \alpha_k} ( \xx_k, \nabla f(\xx_k) )
	\ea
	$$ \\[5pt]
	\ea\\
	\hline
	\ea
	\eeq
	
	Note that when $\psi \equiv 0$, composite step \eqref{CompositeStep}
	coincides with our basic preconditioned gradient step defined earlier:
	$$
	\ba{rcl}
	\text{CompositeStep}_{\mat{H}, \alpha}(\xx, \gg) & = & 
	\text{GradStep}_{\mat{H}, \alpha}(\xx, \gg) 
	\;\; = \;\;
	\xx - \bigl( \mat{H} + \alpha \mat{I}  \bigr)^{-1} \mat{g},
	\ea
	$$ 
	and method~\eqref{CompositeMethod} is identical to \eqref{MainAlgorithm}.

	\subsection{Convergence Analysis}
	
	In this section, we provide the proofs of our main convergence results from Section~\ref{SectionSpectral}.
	We study the more general composite formulation~\eqref{CompositeMethod} of our method,
	which covers the basic case when $\psi \equiv 0$.
	
	Let us consider one step of our method: $\xx^+ = \text{CompositeStep}_{\mat{H}, \alpha}(\xx, \nabla f(\xx))$,
	for some $\alpha > 0$, and establish its key properties.
	The new point $\xx^+$ satisfies the following optimality condition (see, e.g., Theorem 3.1.23 in \cite{nesterov2018lectures}):
	\beq \label{StatCondition}
	\ba{rcl}
	\la \nabla f(\xx) + (\mat{H} + \alpha \mat{I})(\xx^+ - \xx), \yy - \xx^+ \ra + \psi(\yy) & \geq & \psi(\xx^+),
	\qquad \forall \yy \in Q.
	\ea
	\eeq
	In other words, the vector $\psi'(\xx^+) := -\nabla f(\xx) - (\mat{H} + \alpha \mat{I})(\xx^+ - \xx)$ belongs to the subdifferential of $\psi$
	at new point:
	$$
	\ba{rcl}
	\psi'(\xx^+) & \in & \partial \psi(\xx^+).
	\ea
	$$
	We denote
	\beq \label{FSubgrDef}
	\ba{rcl}
	F'(\xx^+) & := & \nabla f(\xx^+) + \psi'(\xx^+),
	\ea
	\eeq
	that is the main object for which we prove the convergence of our method.
	Utilizing positive semi-definitiveness of $\mat{H} = \mat{H}^{\top} \succeq 0$,
	we can can bound the length of our displacement $r \Def \|\xx^+ - \xx\|$.
	
	\BL For any $\ss \in \partial \psi(\xx)$, we have
	\beq \label{RBound}
	\ba{rcl}
	r & \leq & \frac{ \| \nabla f(\xx)  + \ss\|}{\alpha},
	\ea
	\eeq
	and
	\beq \label{HessBound}
	\ba{rcl}
	\la \mat{H}(\xx^+ - \xx), \xx^+ - \xx \ra
	& \leq & r \|  \nabla f(\xx) + \ss\|.
	\ea
	\eeq
	\EL
	\proof
	Indeed, using convexity of $\psi$,
	we obtain, for any $\ss \in \partial \psi(\xx)$
	and our specific $\psi'(\xx^+) \in \partial \psi(\xx^+)$:
	$$
	\ba{rcl}
	0 & \leq & \la \psi'(\xx^+) -  \ss, \xx^+ - \xx \ra \\
	\\
	& = & 
	\la - \nabla f(\xx) - ( \mat{H} + \alpha \mat{I} )(\xx^+ - \xx) - \ss, \xx^+ - \xx \ra
	\\
	\\
	& = &
	\la \nabla f(\xx) + \ss, \xx - \xx^+ \ra - \la (\mat{H} + \alpha \mat{I})(\xx^+ - \xx), \xx^+ - \xx \ra.  
	\ea
	$$ 
	Rearranging the terms and using Cauchy-Schwartz inequality,
	we get
	$$
	\ba{rcl}
	r \|  \nabla f(\xx) + \ss  \| & \geq & 
	\la  \nabla f(\xx) + \ss , \xx - \xx^+ \ra
	\;\; \geq \;\;
	\la \mat{H}(\xx^+ - \xx), \xx^+ - \xx \ra + \alpha r^2.
	\ea
	$$ 
	Taking into account that $\mat{H} \succeq \0$ completes the proof.
	\qed
	
	Therefore, by choosing regularization parameter $\alpha$ appropriately,
	we can control the length of steps for our algorithm.
	When combined with smoothness properties \eqref{LipHessian}, \eqref{SigmaHessian}
	of the objective, we can establish the global progress in terms of the objective function value.
	We denote 
	$$
	\ba{rcl}
	\delta & := & \| \mat{H} - \nabla^2_{\tau} f(\xx) \|.
	\ea
	$$
	
	\BL \label{LemmaCompositeFuncStep}
	Let $\ss \in \partial \psi(\xx)$ be some subgradient of $\psi$ at current point $\xx$
	and
	let $\alpha \geq \sqrt{\frac{L \| \nabla f(\xx) + \ss  \|}{3}} + \sigma_{\tau} + \delta$. Then
	\beq \label{FuncProgress}
	\ba{rcl}
	F(\xx) - F(\xx^+) & \geq & \frac{\alpha}{2} r^2.
	\ea
	\eeq
	\EL
	\proof
	By Lipschitz continuity of the Hessian, we have
	$$
	\ba{rcl}
	F(\xx^+) & = & f(\xx^+) + \psi(\xx^+) \\
	\\
	& \leq & f(\xx) + \la \nabla f(\xx), \xx^+ - \xx \ra
	+ \frac{1}{2} \la \nabla^2 f(\xx)(\xx^+ - \xx), \xx^+ - \xx \ra
	+ \frac{L}{6} r^3 + \psi(\xx^+)\\
	\\
	& \overset{\eqref{RBound}}{\leq} & 
	f(\xx) + \la \nabla f(\xx), \xx^+ - \xx \ra
	+ \frac{1}{2} \la \nabla^2 f(\xx)(\xx^+ - \xx), \xx^+ - \xx \ra
	+ \frac{L \| \nabla f(\xx) + \ss \| }{6 \alpha} r^2 + \psi(\xx^+) \\
	\\
	& \leq & 
	f(\xx) + \la \nabla f(\xx), \xx^+ - \xx \ra
	+ \frac{1}{2} \la \mat{H}(\xx^+ - \xx), \xx^+ - \xx \ra
	+ \frac{1}{2} \bigl( \sigma_{\tau} + \delta
	+ \frac{L \| \nabla f(\xx) + \ss \|}{3\alpha} \bigr) r^2
	+ \psi(\xx^+),
	\ea
	$$
	where $\sigma_{\tau} := \sup_{\xx} \| \nabla^2 f(\xx) - \nabla^2_{\tau} f(\xx) \|$
	as defined in \eqref{SigmaHessian}.
	According to~\eqref{StatCondition}, it holds
	$$
	\ba{rcl}
	\psi(\xx^+) & \leq & \psi(\xx) 
	- \la \nabla f(\xx), \xx^+ - \xx \ra
	- \la \mat{H}(\xx^+ - \xx), \xx^+ - \xx \ra
	- \alpha r^2.
	\ea
	$$
	Hence, substituting this inequality into the previous one, we can continue as follows:
	$$
	\ba{rcl}
	F(\xx^+) & \leq & 
	F(\xx) 
	- \frac{\alpha}{2} r^2 
	- \frac{1}{2} \bigl( \alpha - \sigma_{\tau} - \delta - \frac{L \| \nabla f(\xx) + \ss \|}{3\alpha}  \bigr) r^2
	- \frac{1}{2} \la \mat{H}(\xx^+ - \xx), \xx^+ - \xx \ra  \\
	\\
	& \leq & 
	F(\xx) - \frac{\alpha}{2} r^2  - \frac{1}{2} \bigl(  \alpha - \sigma_{\tau} - \delta -  \frac{L \| \nabla f(\xx) + \ss\|}{3\alpha} \bigr) r^2.
	\ea
	$$
	To prove the result, it suffices to check $\alpha \geq \sigma_{\tau} + \delta + \frac{L \| \nabla f(\xx) + \ss \|}{3\alpha}$,
	which is ensured by our choice. Indeed,
	denoting $B := \sigma_{\tau} + \delta$ and $C := \frac{L \| \nabla f(\xx)  + s\|}{2}$,
	the inequality which we need to ensure for $\alpha > 0$ is
	\beq \label{AlphaCondition}
	\ba{rcl}
	\alpha^2 & \geq & B + \frac{C}{\alpha}
	\quad \Leftrightarrow \quad
	\alpha^2 - \alpha B - C \;\; \geq \;\; 0
	\quad \Leftrightarrow \quad
	\alpha \;\; \geq \;\;
	\frac{B + \sqrt{B^2 + 4C}}{2}.
	\ea
	\eeq
	It is immediate to check that $\alpha := B + \sqrt{C}$
	satisfies \eqref{AlphaCondition}, hence
	it holds for any $\alpha \geq B + \sqrt{C}$, which is our choice in the condition of the lemma.
	\qed

	Now, let us relate the length $r$ of the step  with the norm of vector $F'(\xx^+) := \nabla f(\xx^+) + \psi'(\xx^+)$.
	\BL
		Let $\ss \in \partial \psi(\xx)$ be some subgradient of $\psi$ at current point $\xx$
	and
	let $\alpha \geq \sqrt{\frac{L \| \nabla f(\xx) + \ss  \|}{2}} + \sigma_{\tau} + \delta$.
	Then
	\beq \label{NewGradNorm}
	\ba{rcl}
	r & \geq & \frac{1}{2\alpha}\| F'(\xx^+) \|.
	\ea
	\eeq
	\EL
	\proof
	Using the definition of $\psi'(\xx^+) := -\nabla f(\xx) - (\mat{H} + \alpha \mat{I})(\xx^+ - \xx)$ 
	and Lipschitz continuity of the Hessian, we have
	$$
	\ba{rcl}
	\| F'(\xx^+) \| & = & \| \nabla f(\xx^+) + \psi'(\xx^+) \| \\
	\\
	& = &
	\| \nabla f(\xx^+) - \nabla f(\xx) - (\mat{H} + \alpha \mat{I})(\xx^+ - \xx) \| \\
	\\
	& \leq &
	\| \nabla f(\xx^+) - \nabla f(\xx) - \nabla^2 f(\xx)(\xx^+ - \xx) \|
	+ \bigl( \alpha + \sigma_{\tau} + \delta \bigr) r \\
	\\
	& \overset{\eqref{LipHessian}}{\leq} &
	\frac{L}{2}r^2 + \bigl( \alpha + \sigma_{\tau} + \delta \bigr) r \\
	\\
	& \overset{\eqref{RBound}}{\leq} &
	\bigl(\alpha + \sigma_{\tau} + \delta + \frac{L \| \nabla f(\xx) + \ss \|}{2 \alpha} \bigr) r
	\;\; \leq \;\; 2\alpha r,
	\ea
	$$
	where the last inequality holds due to $\alpha \geq \sigma_{\tau} + \delta + \frac{L \| \nabla f(\xx) + \ss \|}{2 \alpha}$,
	which is ensure by our choice (see the end of the proof of Lemma~\ref{LemmaCompositeFuncStep}).
	\qed
	
	Therefore, combining these two Lemmas together, we obtain the following bound
	for one step of the method.
	\BC
	Let $\alpha \geq \sqrt{\frac{L \| \nabla f(\xx) + \ss \|}{2}} + \sigma_{\tau} + \delta$ for some $\ss \in \partial \psi(\xx)$. Then
	\beq \label{FuncProgress2}
	\ba{rcl}
	F(\xx) - F(\xx^+) & \overset{\eqref{FuncProgress}, \eqref{NewGradNorm}}{\geq} &
	\frac{1}{8 \alpha} \| F'(\xx^+)\|^2. 
	\ea
	\eeq
	\EC
	
	We are ready to prove the global complexity bound for convergence of our method.
	\BT \label{TheoremCompositeLip}
	Let $f \in \F_{\tau}$ be non-convex of grade $\tau$, where $0 \leq \tau \leq n$
	is fixed.
	Let $f$ have a Lipschitz Hessian with constant $L$
	and bounded parameter $\sigma_{\tau}$. 
	Consider iterations $\{ \xx_k \}_{k \geq 0}$
	of algorithm~\eqref{CompositeMethod} with
	$$
	\ba{rcl}
	\alpha_k & = & \sqrt{\frac{L \| F'(\xx_k) \|}{2}} + \sigma_{\tau} + \delta,
	\ea
	$$
	where $F'(\xx_k)$ is defined by \eqref{FSubgrDef} for $k \geq 1$, and
	$F'(\xx_0) := \nabla f(\xx_0) + \ss$, for an arbitrary initial subgradient $\ss \in \partial \psi(\xx_0)$.
	Then, for any $\varepsilon > 0$, it is enough to do
	\beq \label{MainResult}
	\ba{rcl}
	K & = & 
	\Bigl\lceil
	
	8 ( F(\xx_0) - F^{\star} ) \cdot \Bigl(   \sqrt{\frac{L}{2}} \frac{1}{\varepsilon^{3/2}} + \frac{\sigma_{\tau} + \delta}{\varepsilon^2}  \Bigr)
	
	\; + \; 2 \ln \frac{\| F'(\xx_0) \|}{\varepsilon}
	\Bigr\rceil.
	\ea
	\eeq
	iterations to have $\min\limits_{1 \leq i \leq K} \| F'(\xx_i) \| \leq \varepsilon$.
	\ET
	\proof
	Let us fix some $k \geq 0$ and assume that for any $0 \leq i \leq k$ we have $g_i \Def \| F'(\xx_i) \| \geq \varepsilon$.
	
	According to \eqref{FuncProgress2}, we obtain, for $0 \leq i \leq k - 1$:
	\beq \label{AlgOneStep}
	\ba{rcl}
	F(\xx_i) - F(\xx_{i + 1}) & \geq & \frac{1}{8\alpha_k} g_{i + 1}^2
	\;\; = \;\;
	\frac{1}{8}\bigl[  \frac{g_{i + 1}}{g_i}  \bigr]^2
	\cdot \Bigl(  \sqrt{\frac{L}{2}} \frac{1}{g_{i + 1}^{3/2}} + \frac{\sigma_{\tau} + \delta}{g_{i + 1}^2}  \Bigr)^{-1} \\
	\\
	& \geq & 
	\frac{1}{8}\bigl[  \frac{g_{i + 1}}{g_i}  \bigr]^2
	\cdot \Bigl(  \sqrt{\frac{L}{2}} \frac{1}{\varepsilon^{3/2}} + \frac{\sigma_{\tau} + \delta}{\varepsilon^2}  \Bigr)^{-1}.
	\ea
	\eeq
	Denote
	$$
	\ba{rcl}
	c & := &\frac{1}{8}\Bigl(  \sqrt{\frac{L}{2}} \frac{1}{\varepsilon^{3/2}} 
	+ \frac{\sigma_{\tau} + \delta}{\varepsilon^2}  \Bigr)^{-1}.
	\ea
	$$
	Then, telescoping bound \eqref{AlgOneStep} and using the inequality between arithmetic and geometric means, we get
	$$
	\ba{rcl}
	F(\xx_0) - F^{\star} & \geq & F(\xx_0) - F(\xx_{k})
	\;\; \overset{\eqref{AlgOneStep}}{\geq} \;\;
	c
	\sum\limits_{i = 0}^{k - 1} \bigl[  \frac{g_{i + 1}}{g_i} \bigr]^2 
	\;\; \geq \;\;
	ck
	\Bigl[\prod\limits_{i = 0}^{k - 1}  \frac{g_{i + 1}}{g_i} \Bigr]^{2 / k} \\
	\\
	& = & 
	ck \Bigl[ \frac{g_k}{g_0}  \Bigr]^{2 / k}
	\;\; \geq \;\;
	ck \Bigl[ \frac{\varepsilon}{g_0}  \Bigr]^{2 / k}
	\;\; = \;\;
	ck \exp\Bigl[ -\frac{2}{k} \ln \frac{g_0}{\varepsilon}  \Bigr] \\
	\\
	& \geq &
	ck \Bigl[ 1 - \frac{2}{k} \ln \frac{g_0}{\varepsilon}  \Bigr]
	\;\; = \;\;
	ck - 2c \ln \frac{g_0}{\varepsilon}.
	\ea
	$$
	Hence, we obtain
	$$
	\ba{rcl}
	k & \leq & \frac{F(\xx_0) - F^{\star}}{c} + 2 \ln \frac{g_0}{\varepsilon}.
	\ea
	$$
	Substituting the value of $c$ completes the proof.
	\qed

	\subsection{Proof of Theorem~\ref{TheoremLip}}
	\proof
	It follows immediately from Theorem~\ref{TheoremCompositeLip}
	by substituting the non-composite case $\psi \equiv 0$. 
	\qed

	\section{Proofs for Section~\ref{SectionCut}}
	\label{SectionAppendixCut}
	
	The results of this section are applied to a basic unconstrained minimization
	problem~\eqref{MainProblem}.
	Let us consider one iteration $\xx \mapsto \xx^+$
	of our method, which satisfies the following stationary condition
	\beq \label{UnStationary}
	\ba{rcl}
	\nabla f(\xx) + \mat{H}(\xx^+ - \xx) + \alpha \mat{I} & = & \0,
	\ea
	\eeq
	where regularization parameter is chosen as 
	\beq \label{UnAlphaChoice}
	\ba{rcl}
	\alpha & := & \alpha^{\star} + \eta,
	\ea
	\eeq
	with
	\beq \label{AlphaStar}
	\ba{rcl}
	\alpha^{\star} & := &
	\argmax\limits_{\alpha > 0}
	\Bigl[ 
	- \frac{1}{2}\la (\mat{H} + (\alpha + \eta) \mat{I} )^{-1} \nabla f(\xx), \nabla f(\xx) \ra
	- \frac{2\alpha^3}{3L}
	\Bigr],	
	\ea
	\eeq
	and $\eta := \sigma_{\tau}^+ + \delta+ \delta_{-} \geq 0$ is the balancing term to control the errors.
	Considering the decomposition of the Hessian onto positive and negative components:
	\beq \label{HessDecompose}
	\ba{rcl}
	\nabla^2 f(\xx) & \equiv & \nabla^2_+ f(\xx) - \nabla^2_{-} f(\xx), 
	\qquad 
	\nabla^2_+ f(\xx), \,
	\nabla^2_- f(\xx)  \;\; \succeq \;\; \mat{0},
	\ea
	\eeq
	we have that
	\beq \label{ErrorBounds}
	\ba{rcl}
	\| \nabla^2_+ f(\xx) - \nabla^2_{\tau} f(\xx) \| & \leq & \sigma_{\tau}^+ \\
	\\
	\| \mat{H} - \nabla^2_{\tau} f(\xx) \| & \leq & \delta, \\
	\\ 
	\| \nabla^2_{-} f(\xx) \mat{P} \| & \leq & \delta_{-},
	\ea
	\eeq
	where $\mat{P}$ is a given projector onto image of $\mat{H}$, that satisfies
	\beq \label{PProperty}
	\ba{rcl}
	\mat{H} \mat{P} & = & \mat{H}.
	\ea
	\eeq
	
	First, let us provide another description of regularization parameter $\alpha^{\star}$,
	that is more suitable for our analysis.
	First-order stationary condition for \eqref{AlphaStar} gives
	$$
	\ba{rcl}
	\frac{1}{2} \la (\mat{H} + \alpha \mat{I} )^{-1} \nabla f(\xx),
	(\mat{H} + \alpha \mat{I} )^{-1} \nabla f(\xx) \ra
	- \frac{2 (\alpha^{\star})^2}{L} \;\; & = & 0,
	\ea
	$$
	which is equivalent to
	\beq \label{AlphaStarFormula}
	\ba{rcl}
	\alpha^{\star} & = & \frac{L}{2}\| (\mat{H} + \alpha \mat{I})^{-1}  \nabla f(\xx) \|
	\;\; \overset{\eqref{UnStationary}}{=} \;\;
	\frac{L}{2}\|\xx^+ - \xx\|
	\;\; \equiv \;\;
	\frac{Lr}{2},
	\ea
	\eeq
	where we use $r := \|\xx^+ - \xx\|$ as previously.
	Therefore, from \eqref{AlphaStarFormula}
	we see that $\alpha^{\star}$ plays the role
	of the Cubic Regularization~\cite{nesterov2006cubic}
	of our model.
	
	Let us establish the main inequalities
	on the progress of each step.
	
	\BL \label{LemmaUnFunc}
	For the functional residual, it holds
	\beq \label{UnFuncRes}
	\ba{rcl}
	f(\xx) - f(\xx^+) & \geq & \frac{L}{3} r^3 + \frac{\eta}{2}r^2
	\ea
	\eeq
	\EL
	\proof
	Indeed, using Lipschitz continuity of the Hessian 
	and definition of our step, we obtain
	$$
	\ba{rcl}
	f(\xx^+) & \overset{\eqref{LipHessian}, \eqref{HessDecompose}}{\leq} &
	f(\xx) + \la \nabla f(\xx), \xx^+ - \xx \ra
	+ \frac{1}{2}\la \nabla^2_+ f(\xx) (\xx^+ - \xx), \xx^+ - \xx \ra
	- \frac{1}{2}\la \nabla^2_{-} f(\xx) (\xx^+ - \xx), \xx^+ - \xx \ra
	+ \frac{L}{6} r^3 \\
	\\
	& \leq & 
	f(\xx)
	+ \la \nabla f(\xx), \xx^+ - \xx \ra
	+ \frac{1}{2}\la \nabla^2_+ f(\xx) (\xx^+ - \xx), \xx^+ - \xx \ra
	+ \frac{L}{6} r^3 \\
	\\
	& \overset{\eqref{StatCondition}}{=} &
	f(\xx) - \la \mat{H}(\xx^+ - \xx), \xx^+ - \xx \ra
	- \alpha r^2
	+ \frac{1}{2}\la \nabla^2_+ f(\xx) (\xx^+ - \xx), \xx^+ - \xx \ra
	+ \frac{L}{6} r^3 \\
	\\
	& \leq &
	- \alpha r^2 + \frac{\sigma_{\tau}^+ + \delta}{2} r^2 + \frac{L}{6} r^3.
	\ea
	$$
	Hence, rearranging the terms and using the definition of $\alpha$, we obtain
	$$
	\ba{rcl}
	f(\xx) - f(\xx^+) & \geq & 
	\alpha r^2 - \frac{\sigma_{\tau}^+ + \delta}{2}r^2 - \frac{L}{6} r^3 \\
	\\
	& \overset{\eqref{UnAlphaChoice}, \eqref{AlphaStarFormula}}{=} & 
	\frac{L}{2} r^3 + \eta r^2 - \frac{\sigma_{\tau}^+ + \delta}{2} r^2 - \frac{L}{6} r^3 \\
	\\
	& \geq &  
	\frac{L}{3} r^3 + \frac{\eta}{2}r^2,
	\ea
	$$
	which is the required bound.
	\qed

	\BL \label{LemmaUnGrad}
	Let $\mat{y}^+ := \mat{x} + \mat{P}(\xx^+ - \xx)$. Then,
	we can relate the gradient at $\yy^+$ and the length $r$ of the step, as follows
	\beq \label{UnGradBound}
	\ba{rcl}
	\| \nabla f(\yy^+) \| & \leq & 2\eta r + L r^2.
	\ea
	\eeq
	\EL
	\proof
	By Lipschitz continuity of the Hessian, we have
	$$
	\ba{rcl}
	\frac{L}{2}r^2 & \geq & 
	\frac{L}{2}\| \yy^+ - \xx \|^2 \\
	\\
	& \overset{\eqref{LipHessian}}{\geq} & 
	\| \nabla f(\yy^+) - \nabla f(\xx) - \nabla^2 f(\xx) (\yy^+ - \xx) \| \\
	\\
	& = & 
	\| \nabla f(\yy^+) - \nabla f(\xx) - \nabla^2 f(\xx) \mat{P} (\xx^+ - \xx) \| \\
	\\
	& \overset{\eqref{StatCondition}}{=} &
	\| \nabla f(\yy^+) + (\mat{H} - \nabla^2 f(\xx)) \mat{P}(\xx^+ - \xx) + \alpha (\xx^+ - \xx) \| \\
	\\
	& \geq &
	\| \nabla f(\yy^+) \| - \| (\mat{H} - \nabla^2 f(\xx)) \mat{P}(\xx^+ - \xx) \| - \alpha r.
	\ea
	$$
	Hence, rearranging the terms, we obtain
	$$
	\ba{rcl}
	\| \nabla f(\yy^+) \| & \leq & \| (\mat{H} - \nabla^2 f(\xx)) \mat{P}(\xx^+ - \xx) \| + \alpha r + \frac{L}{2} r^2 \\
	\\
	& \overset{\eqref{HessDecompose}}{=} &
	\| (\mat{H} - \nabla^2_+ f(\xx) + \nabla^2_{-} f(\xx)) \mat{P}(\xx^+ - \xx) \|  + \alpha r + \frac{L}{2} r^2 \\
	\\
	& \leq & 2\eta r + Lr^2,
	\ea
	$$
	which is the required bound.
	\qed
	
	Combining inequalities \eqref{UnFuncRes} and \eqref{UnGradBound} together, we obtain the following bound
	on the progress.
	\BC
	For one step of the method,
	it holds:
	\beq \label{UnStepProgress}
	\ba{rcl}
	f(\xx) - f(\xx^+) & \geq & 
	\min\Bigl\{  \frac{1}{6\sqrt{2L}}\| \nabla f(\yy^+) \|^{3/2}, \, 
	\frac{1}{32 \eta} \| \nabla f(\yy^+)\|^2    \Bigr\}.
	\ea
	\eeq
	\EC

	\subsection{Proof of Theorem~\ref{TheoremCut}}
	\label{SectionAppendixProofCut}
	
	\proof
	
	Let us fix some $k \geq 1$ and assume that for any $1 \leq i \leq k$
	we have $\| \nabla f(\yy_i) \| \geq \varepsilon$.
	
	According to \eqref{UnStepProgress}, we obtain, for $0 \leq i \leq k - 1$:
	$$
	\ba{rcl}
	f(\xx_i) - f(\xx_{i + 1}) & \geq & 
	\min\Bigl\{  \frac{1}{6\sqrt{2L}}\| \nabla f(\yy_{i + 1}) \|^{3/2}, \, 
	\frac{1}{32 \eta} \| \nabla f(\yy_{i + 1})\|^2    \Bigr\} \\
	\\
	& \geq & 
	\min\Bigl\{  \frac{1}{6\sqrt{2L}}\varepsilon^{3/2}, \, 
	\frac{1}{32 \eta} \varepsilon^2    \Bigr\}.
	\ea
	$$
	Telescoping this bound for the first $k$ iterations, we get
	$$
	\ba{rcl}
	f(\xx_0) - f^{\star} & \geq &
	f(\xx_0) - f(\xx_k) \;\; \geq \;\;
	k 	\min\Bigl\{  \frac{1}{6\sqrt{2L}}\varepsilon^{3/2}, \, 
	\frac{1}{32 \eta} \varepsilon^2    \Bigr\},
	\ea
	$$
	which leads to the required complexity.
	\qed
	
	\section{Proofs for Section~\ref{SectionConvex}}
	\label{SectionAppendixConvex}
	
	In this section, we provide a general analysis for
	the \textit{composite} version of our method~\eqref{CompositeMethod},
	when the target objective is convex: $f \in \F_n$
	and quasi-Self-Concordant~\eqref{MquasiSC} with parameter $M > 0$.

	Let us establish the progress for one step of our algorithm 
	$\xx^+ = \text{CompositeStep}_{\mat{H},\alpha}(\xx, \nabla f(\xx))$
	under this refined smoothness condition.
	
	\BL
	Let $\alpha \geq M \| F'(\xx) \| + \sigma_{\tau} + \delta$. Then
	\beq \label{QSCOneStep}
	\ba{rcl}
	\la F'(\xx^+), \xx - \xx^+ \ra & \geq & \frac{1}{2\alpha} \| F'(\xx^+) \|^2.
	\ea
	\eeq
	\EL
	\proof
	By using basic properties of quasi-Self-Concordant functions (see Lemma 2.7 in \cite{doikov2023minimizing}), we have, 
	for any two points $\xx, \xx^+ \in \R^n$:
	$$
	\ba{rcl}
	\| \nabla f(\xx^+) -  \nabla f(\xx) - \nabla^2 f(\xx)(\xx^+ - \xx) \| & \leq &
	\frac{M}{2} \la \nabla^2 f(\xx)(\xx^+ - \xx), \xx^+ - \xx \ra \cdot \varphi( Mr ),
	\ea
	$$
	where $\varphi(t) \Def \frac{e^t - t - 1}{t^2}$ is a convex monotone univariate function.
	Let us show how we can control the right hand side.
	From~\eqref{RBound}, we have
	$$
	\ba{rcl}
	\varphi( Mr ) & \leq & \varphi( \frac{M \| F'(\xx) \|}{\alpha} )
	\;\; \overset{(*)}{\leq} \;\; 1,
	\ea
	$$
	where in $(*)$ we use our choice of the regularization coefficient: $\alpha \geq M \| F'(\xx) \|$.
	
	For the Hessian, we can use
	$$
	\ba{rcl}
	\la \nabla^2 f(\xx)(\xx^+ - \xx), \xx^+ - \xx \ra
	& \overset{\eqref{SigmaHessian}}{\leq} &
	\sigma_{\tau} r^2 + \la \mat{H}(\xx^+ - \xx), \xx^+ - \xx \ra
	\;\; \overset{\eqref{HessBound}}{\leq} \;\;
	\sigma_{\tau} r^2 + r \| F'(\xx) \| \\
	\\
	& \overset{\eqref{RBound}}{\leq} & 
	(\frac{\sigma_{\tau}}{\alpha} + 1 )\cdot  r \| F'(\xx) \|
	\;\; \overset{(**)}{\leq} \;\; 2 r \| F'(\xx) \|,
	\ea
	$$
	where we used in $(**)$ that $\alpha \geq \sigma_{\tau}$.
	
	Therefore, we obtain
	\beq \label{GradLinAppxBound}
	\ba{rcl}
	\| \nabla f(\xx^+) - \nabla f(\xx) - \nabla^2 f(\xx)(\xx^+ - \xx) \|
	& \leq & 
	M r \| F'(\xx) \|.
	\ea
	\eeq
	Thus, using the definition of new subgradient $F'(\xx^+)$ from the stationary condition~\eqref{StatCondition}, we have
	$$
	\ba{rcl}
	\| F'(\xx^+) + \alpha (\xx^+ - \xx) \|
	& = &
	\| \nabla f(\xx^+) - \nabla f(\xx) - \mat{H}(\xx^+ - \xx) \| \\
	\\
	& \leq & 
	\| \nabla f(\xx^+) - \nabla f(\xx) - \nabla^2 f(\xx)(\xx^+ - \xx) \| \\
	\\
	& & \;\; + \;\; r \| \mat{H} - \nabla^2 f(\xx) \| \\
	\\
	& \overset{\eqref{GradLinAppxBound}}{\leq} &
	r(\sigma_{\tau} + \delta + M \| F'(\xx)\|).
	\ea
	$$
	Taking square of both sides and rearranging the terms, we obtain
	$$
	\ba{rcl}
	\la F'(\xx^+), \xx - \xx^+ \ra 
	& \geq & 
	\frac{1}{2 \alpha} \| F'(\xx^+) \|^2
	+ \frac{r^2}{2 \alpha} (\alpha^2 - (\sigma_{\tau} + \delta + M \| \nabla f(\xx) \|)^2).
	\ea
	$$
	Taking into account or choice of $\alpha$ completes the proof.
	\qed

	We denote by $D$ the diameter of the initial sublevel set which we assume to be bounded:
	$$
	\ba{rcl}
	D & := &
	\sup\limits_{\xx} \Bigl\{
	\| \xx - \xx^{\star} \| \; : \; F(\xx) \leq F(\xx_0)
	\Bigr\} 
	\;\; < \;\; +\infty.
	\ea
	$$
	We prove the following result. 
	
	\BT \label{TheoremCompositeConvex}
	Let $f \in \F_n$ be convex and quasi-Self-Concordant wtih constant $M$.
	Consider iterations $\{ \xx_k \}_{k \geq 0}$ of algorithm~\eqref{CompositeMethod}
	with
	$$
	\ba{rcl}
	\alpha_k & = & M \| F'(\xx_k) \| + \sigma_\tau + \delta,
	\ea
	$$
	where $F'(\xx_k)$ is defined by \eqref{FSubgrDef} for $k \geq 1$, and
	$F'(\xx_0) := \nabla f(\xx_0) + \ss$, for an arbitrary initial subgradient $\ss \in \partial \psi(\xx_0)$.
	Then, for any $\varepsilon > 0$, it is enough to do
	\beq \label{AppendixConvexComplexity}
	\ba{rcl}
	K & = &
	4 \Bigl\lceil
	\Bigl( MD + \frac{\sigma_{\tau} + \delta}{2\mu}  \Bigr)
	\ln \frac{F(\xx_0) - F^{\star}}{\varepsilon}
	\; + \; \ln \frac{\| F'(\xx_0)\| D }{\varepsilon}
	\Bigr\rceil
	\ea
	\eeq
	steps to ensure $F(\xx_k) - F^{\star} \leq \varepsilon$.
	\ET
	\begin{proof}
		Let us prove the following rate of convergence, for the iterations of our method, 
		\beq \label{MainRate}
		\ba{rcl}
		F(\xx_k) - F^{\star} & \leq & 
		\exp\Bigl( -\frac{k}{4} \Bigl[ \frac{\sigma_{\tau} + \delta}{2\mu} + MD   \Bigr]^{-1}  \Bigr)
		\bigl( F(\xx_0) - F^{\star}   \bigr)
		\; + \; \exp\Bigl( -\frac{k}{4}  \Bigr) \cdot \| F'(\xx_0) \| D,
		\ea
		\eeq
		which immediately leads to the complexity bound~\eqref{AppendixConvexComplexity}.
		
		We denote $g_k = \| F'(\xx_k)\|$.
		Then, by convexity, we have for every iteration $k \geq 0$:
		\beq \label{FuncProgress3}
		\ba{rcl}
		F(\xx_k) - F(\xx_{k + 1}) & \geq & 
		\la F'(\xx_{k + 1}), \xx_k - \xx_{k + 1} \ra \\
		\\
		& \overset{\eqref{QSCOneStep}}{\geq} & 
		\frac{1}{2 ( \sigma_{\tau} + \delta + M g_k  )} \bigl( \frac{g_{k + 1}}{g_k}  \bigr)^2 g_k^2 \\
		\\
		& = & 
		\frac{1}{2} \bigl( \frac{g_{k + 1}}{g_k}  \bigr)^2 \cdot \Bigl(  \frac{\sigma_{\tau} + \delta}{g_k^2} + \frac{M}{g_k}  \Bigr)^{-1}.
		\ea
		\eeq
		Denoting the functional residual by $F_k \Def F(\xx_k) - F^{\star}$,
		we have by convexity and strong convexity:
		$$
		\ba{rcl}
		g_k & \geq & \frac{F_k}{D}, \qquad g_k^2 \;\; \geq \;\; 2\mu F_k.
		\ea
		$$
		Substituting these bounds into \eqref{FuncProgress3}, we obtain
		\beq \label{FuncProgress4}
		\ba{rcl}
		F_k - F_{k + 1} & \geq & 
		\frac{1}{2} \bigl( \frac{g_{k + 1}}{g_k}  \bigr)^2 
		\cdot \Bigl(  
		\frac{\sigma_{\tau} + \delta}{2\mu F_k} + \frac{MD}{F_k} 
		\Bigr)^{-1}
		\;\; = \;\; \frac{q}{2} \bigl( \frac{g_{k + 1}}{g_k}  \bigr)^2 
		F_k,
		\ea
		\eeq
		where condition number $q$ is defined by
		$$
		\ba{rcl}
		q &:=& \Bigl( \frac{\sigma_{\tau} + \delta}{2\mu} + MD \Bigr)^{-1}.
		\ea
		$$
		To show the desired rate, we use concavity of logarithm,
		$$
		\ba{rcl}
		\ln \frac{a}{b} & = & \ln a - \ln b \;\; \geq \;\; \frac{1}{a}(a - b), \qquad \forall a, b > 0,
		\ea
		$$
		and conclude that
		$$
		\ba{rcl}
		\ln \frac{F_k}{F_{k + 1}} & \geq & \frac{1}{F_k}(F_k - F_{k + 1})
		\;\; \overset{\eqref{FuncProgress4}}{\geq} \;\;
		\frac{q}{2} \bigl(  \frac{g_{k + 1}}{g_k} \bigr)^2.
		\ea
		$$
		Telescoping this bound and using inequality between arithmetic and geometric means, we get
		$$
		\ba{rcl}
		\ln \frac{F_0}{F_k} & \geq & 
		\frac{q}{2} \sum\limits_{i = 0}^{k - 1} \bigl[  \frac{g_{i + 1}}{g_i} \bigr]^2 
		\;\; \geq \;\;
		\frac{kq}{2} \Bigl[  \prod\limits_{i = 0}^{k - 1} \frac{g_{i + 1}}{g_i} \Bigr]^{2/k}
		\;\; = \;\;
		\frac{kq}{2} \bigl[ \frac{g_k}{g_0} \bigr]^{2/k} \\
		\\
		& = &
		\frac{kq}{2} \exp \bigl[ \frac{2}{k} \ln \frac{g_k}{g_0} \bigr]
		\;\; \geq \;\;
		\frac{kq}{2}  \Bigl( 1 + \frac{2}{k} \ln \frac{g_k}{g_0}  \Bigr)
		\;\; \geq \;\;
		\frac{kq}{2}  \Bigl( 1 + \frac{2}{k} \ln \frac{F_k}{g_0 D}  \Bigr).
		\ea
		$$
		The last inequality leads to the required bound \eqref{MainRate}.
	\end{proof}

	\subsection{Proof of Theorem~\ref{TheoremConvex}}
	\proof
	It follows immediately from Theorem~\ref{TheoremCompositeConvex}
	by substituting the non-composite case $\psi \equiv 0$. 
	\qed

\end{document}